\documentclass[11pt]{article}

\usepackage[truedimen,margin=25truemm]{geometry}

\usepackage[dvipdfmx]{graphicx}

\usepackage{additional_def}

\title{Fast online low-rank tensor subspace tracking by\\ CP decomposition using recursive least squares\\from incomplete observations}
\date{\today}

\author{Hiroyuki Kasai\thanks{Graduate School of Informatics and Engineering, The University of Electro-Communications, Tokyo, Japan ({\tt kasai@is.uec.ac.jp}).} }

\begin{document}

\maketitle

\begin{abstract}
We consider the problem of online subspace tracking of a partially observed high-dimensional data stream corrupted by noise, where we assume that the data lie in a low-dimensional linear subspace. This problem is cast as an online low-rank tensor completion problem. We propose a novel online tensor subspace tracking algorithm based on the CANDECOMP/PARAFAC (CP) decomposition, dubbed OnLine Low-rank Subspace tracking  by TEnsor CP Decomposition (OLSTEC). The proposed algorithm especially addresses the case in which the subspace of interest is dynamically time-varying. To this end, we build up our proposed algorithm exploiting the recursive least squares (RLS), which is the second-order gradient algorithm. Numerical evaluations on synthetic datasets and real-world datasets such as communication network traffic, environmental data, and surveillance videos, show that the proposed OLSTEC algorithm outperforms state-of-the-art online algorithms in terms of the convergence rate per iteration.
\end{abstract}
\section{Introduction}
\label{sec:intro}

The analysis of {\it big data} characterized by a huge volume of massive data is at the very core of recent machine  learning, signal processing, and statistical learning. 
The data have a naturally {\it multi-dimensional} structure. Then they are naturally represented by a multi-dimensional array matrix, namely, a {\it tensor}. 
When the data are high-dimensional data  corrupted by noise, it is very challenging to reveal its underlying latent structure, such as to obtain meaningful information, to impute a missing value, to remove the noise, or to predict some future behaviors of data of interest. 
For this purpose, one typical but promising approach exploits the structural assumption that the data of interest have {\it low-dimensional subspace}, i.e., {\it low-rank}, in every dimension. Many data analysis tasks are achieved efficiently by considering {\it singular value decomposition} (SVD), which reveals the latent subspace of the data.
However, when the data have missing elements caused by, for example, system error, or communication error, SVD cannot be applied directly. To address this shortcoming, low-rank {\it tensor completion} has been studied intensively in recent years. 
A \emph{convex relaxation} \cite{Liu_IEEETransPAMI_2013_s, Tomioka_Latent_2011_s, Signoretto_MachineLearning_2014_s} approach, which is a popular method, estimates the subspace by minimizing the sum of the nuclear norms of the unfolding matrices of the tensor of interest. This approach extends the successful results in the matrix completion problem \cite{Candes_FoundCompuMath_2009} accompanied with theoretical performance guarantees. 
However, because of the high computation cost necessary for the SVD calculation of big matrices every iteration, its scalability is limited on very large-scale data. Instead, a \emph{fixed-rank} non-convex approach with tensor decomposition \cite{Filipovi_MultiSysSigPro_2013_s, Kressner_BIT_2014_s} has gained great attentions recently because of superior performance in practice irrespective of introduction of local minima. This performance also derives from the success of matrix cases \cite{Boumal_NIPS_2011,Mishra_SIAMOpt_2013,Ngo_NIPS_2012_s}. This paper follows the same line as that of the latter approach. 

%
When the data are acquired sequentially from time to time, it is more challenging because of the need for {\it online-based} analysis without storing all of the past data as well as without reliance on the {\it batch-based} process. From this perspective, the batch-based SVD approach is inefficient. It cannot be applied for real-time processing. 
For this problem, {\it online subspace tracking} plays a fundamentally important role in various data analyses to avoid expensive repetitive computations and high memory/storage consumption. 

This present paper particularly addresses two special but realistic situations that arise in the online subspace tracking in practical applications. 
First, (i) considering the {\it time-varying} dynamic nature of real-world streaming data, because there might not exist a {\it unique} and {\it stationary} subspace over time, we are often required to update such a time-varying subspace from moment to moment without sweeping the data in multiple times. Despite allowing moderate accuracy of subspace estimation, this update makes existing batch-based algorithms useless. In fact, as experiments described  later in the paper reveal, such a batch-based approach does not work well under the situation where a stationary subspace does not exist. 
%
Furthermore, (ii) considering the situation and applications where the {\it computational speed} is much faster than the {\it data acquiring speed}, we prefer the algorithm with faster {\it convergence rate} in terms of iteration rather than that with faster computational speed.
For all of these reasons, 
we particularly address the recursive least squares (RLS) algorithm. Although the RLS does not give higher precision from the viewpoint of the optimization theory \cite{Haykin_AdaptiveFilterBook_2002}, it fits the dynamic situation as considered herein because it achieves much faster convergence rate per iteration as a result of the second-order optimization feature. 

This paper presents a new online tensor tracking algorithm, dubbed OnLine Low-rank Subspace tracking  by TEnsor CP Decomposition (OLSTEC), for the partially observed high-dimensional data stream corrupted by noise. We specifically examine the fixed-rank tensor completion algorithm with the second-order gradient descent based on the CP decomposition exploiting the RLS. The advantage of the proposed algorithm, OLSTEC, is quite robust for dynamically time-varying subspace, which often arises in practical applications. This engenders faster update of sudden change of subspaces of interest. This capability is revealed in the numerical experiments conducted with several benchmarks at the end of this paper. 

The remainder of paper is organized as follows. Section 2 introduces some preliminaries followed by related work. Section 3 formulates the problem of online subspace tracking, and proposes details of the new algorithm: OLSTEC. Section 4 presents a description of two extensions of the proposed algorithm, and Section 5 provides computational complexity and memory consumption analysis. Numerical evaluations are performed in Section 6 on synthetic datasets and real-world datasets such as communication network traffic, environmental data and surveillance video. Then, we conclude in Section 7. This paper is an extended version of a short conference paper \cite{Kasai_IEEEICASSP_2016}.

\section{Preliminaries and related work}

\subsection{Preliminaries}
\label{sec:preliminaries}

This subsection first summarizes the notations used in the remainder of this paper. It then briefly introduces the CANDECOMP/PARAFAC (CP) decomposition and the RLS algorithm, which are the basic techniques of the proposed algorithm. 

\subsubsection{Notations}
\label{sec:notations}
We denote scalars by lower-case letters $(a, b, c, \ldots)$, vectors as bold lower-case letters $(\vec{a}, \vec{b}, \vec{c}, \ldots)$, and matrices as bold-face capitals $(\mat{A}, \mat{B}, \mat{C}, \ldots)$. An element at $(i,j)$ of a matrix \mat{A} is represented as $\mat{A}_{i,j}$. 
We use $(\mat{A}[t])_{i,j}$ or $(\mat{BC})_{i,j}$ with a parenthesis if $\mat{A}$ has additional index such as $\mat{A}[t]$ or $\mat{A}$ is a matrix product such as $\mat{A}=\mat{BC}$. 
The $i$-th row vector and $j$-th column of $\mat{A}$ are represented as $\mat{A}_{i,:}$ and $\mat{A}_{:,j}$, respectively. 
It is noteworthy that the transposed column vector of $i$-th row vector $\mat{A}_{i,:}$ is specially denoted as $\vec{a}^i$ with superscript to express a row vector explicitly, i.e., a horizontal vector. $\mat{A}_{i,p:q}$ represents $(\mat{A}_{i,p}, \ldots, \mat{A}_{i,q}) \in \mathbb{R}^{1 \times (q-p+1)}$. 
$\mat{I}_p$ is an identity matrix size of $p\times p$. 
We designate a multidimensional or multi-{\it way} (also called {\it order} or {\it mode}) array as a {\it tensor}, which is denoted by $(\mathbfcal{A}, \mathbfcal{B}, \mathbfcal{C}, \ldots)$. Similarly, an element at $(i,j,k)$ of a third-order tensor $\mathbfcal{A}$ is expressed as $\mathbfcal{A}_{i,j,k}$.
Tensor {\it slice} matrices are defined as two-dimensional matrices of a tensor, defined by fixing all but two indices. 
For example, a {\it horizontal slice} and a {\it frontal slices} of a third-order tensor $\mathbfcal{A}$ are denoted, respectively, as $\mathbfcal{A}_{i,:,:}$ and $\mathbfcal{A}_{:,:,k}$. Also, $\mathbfcal{A}_{:,:,k}$ is use heavily in this article. Therefore, it is simply expressed as $\mat{A}_k$ using the bold-face capital font and a single subscript to  represent its matrix form explicitly. Finally, $\vec{a}[t]$ and $\mat{A}[t]$ with the {\it square bracket} represent the computed $\vec{a}$ and $\mat{A}$ after performing $t$-times updates (iterations) in the online-based subspace tracking algorithm.
The notation {\rm diag}(\vec{a}), where \vec{a} is a vector, stands for the diagonal matrix with $\{\vec{a}_i\}$ as diagonal elements.
We follow the tensor notation of the review article  \cite{Kolda_SIAMReview_2009} throughout our article and refer to it for additional details.
The symbol $\circledast$ denotes the Hadamard Product, which is the element-wise product.

\subsubsection{CANDECOMP/PARAFAC (CP) decomposition}
\label{sec:CPDEC}
\begin{figure}[htbp]
	\begin{center}
	\includegraphics{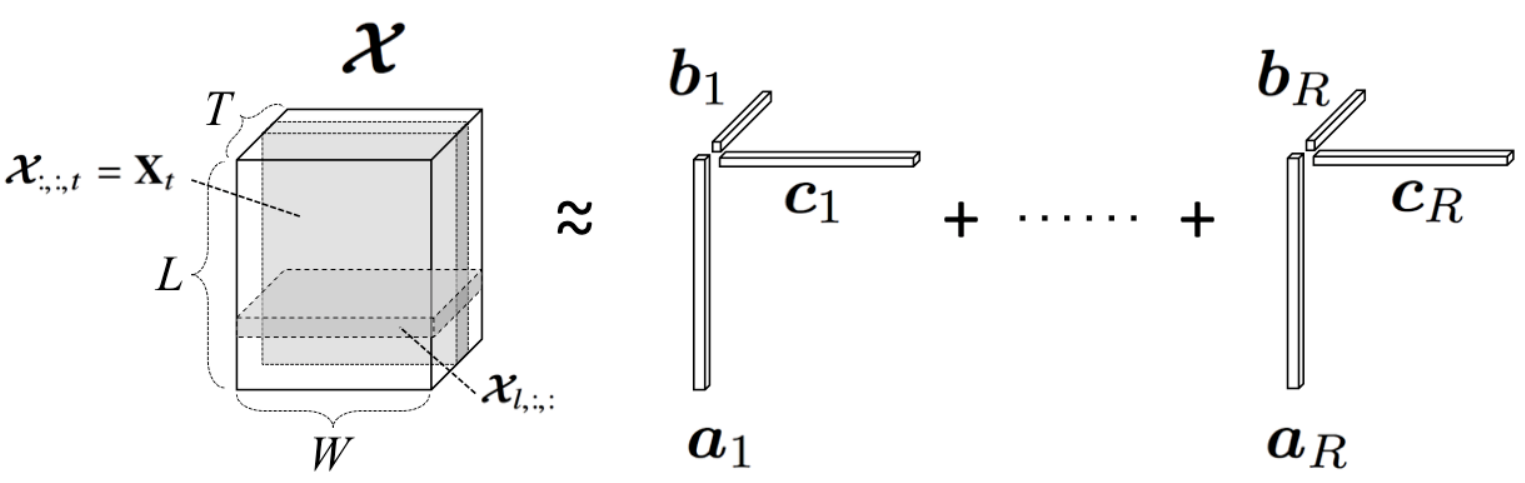}
	\caption{CANDECOMP/PARAFAC (CP) tensor decomposition.}
	\label{Append_Fig:BasicConcept}
	\end{center}
\end{figure}

The CANDECOMP/PARAFAC (CP) decomposition decomposes a tensor into a sum of component rank-one tensors \cite{Kolda_SIAMReview_2009}. Figure \ref{Append_Fig:BasicConcept} presents rank-one tensor decomposition of the Candecomp/PARAFAC decomposition. Letting $\mathbfcal{X}$ be a third-order tensor of size $L\times W \times T$, and assuming its {\it rank} is $R$, we approximate $\mathbfcal{X}$ as $\mathbfcal{X}  \approx \mat{A}{\rm diag}(\vec{b}^t)\mat{C}^T = \sum_{r=1}^R \vec{a}_r \circ \vec{c}_r  \circ \vec{b}_r=\sum_{r=1}^R \vec{b}^t(r)\vec{a}_ r \vec{c}_r^T$, where $\vec{a}_r \in \mathbb{R}^L$, $\vec{b}_r \in \mathbb{R}^W$, and $\vec{c}_r \in \mathbb{R}^T$. The symbol ``$\circ$'' represents the vector outer product. The {\it factor matrices} refer to the combination of the vectors from the rank-one components, i.e., \mat{A} = $[\vec{a}_1: \vec{a}_2: \cdots : \vec{a}_R] \in \mathbb{R}^{L \times R}$ and likewise for $\mat{B} \in \mathbb{R}^{W \times R}$ and $\mat{C }\in \mathbb{R}^{T \times R}$. It must be emphasized that \mat{A}, \mat{B} and \mat{C} can also be represented by {\it row vectors}, i.e., {\it horizontal vectors}, as, 
$\mat{A} = [(\vec{a}^1)^T: \cdots: (\vec{a}^L)^T]^T$, 
$\mat{B} = [(\vec{b}^1)^T: \cdots : (\vec{b}^T)^T]^T$, and
$\mat{C} = [(\vec{c}^1)^T : \cdots : (\vec{c}^W)^T]^T$, where $\{\vec{a}^l, \vec{b}^t,\vec{c}^w\} \in \mathbb{R}^{R}$. 
%
%

\subsubsection{Recursive least squares (RLS)}
\label{Sec:RLS}

\textcolor{black}{
In a least-squares (LS) problem, unknown parameters of a linear model are calculated by minimizing the sum of the squares of the difference between the computed values and the actually observed. To optimize such a least-squares problem, we have a closed form solution. When the interest is in a real-time calculation, it is computationally more efficient if we update the estimates recursively as new data becomes available online. For this purpose, the recursive least squares (RLS) algorithm is a popular algorithms, which is used in adaptive control, adaptive filtering, and system identification \cite{Haykin_AdaptiveFilterBook_2002}. RLS offers a superior convergence rate especially for highly correlated input signals, which is regarded as optimal in practice. This advantage has a price: an increase in the computational complexity. 
Actually, RLS incorporates the history of errors of a considered system into the calculation of the present error compensation. The primary topic of investigation was 
{\it forgetting parameter}, $\lambda$. The forgetting parameter decides how exponentially less important the history of errors is. Although $\lambda=1$ gives the same weights for all the history, the values of $\lambda<1$ bring an exponential decrease in the importance of the history. Finally, it is noteworthy that when implemented in a finite precision environment, the RLS algorithm can suddenly become unstable.  Furthermore, divergence comes to present a difficulty.}

\subsection{Related work}
\label{sec:related_work}

This subsection details general online-based subspace learning methods into which our approach falls. They have been studied actively in the machine learning field recently, and are applicable to noisy, high-dimensional, and incomplete measurements.

Representative research of the {\it matrix-based} online algorithm is the projection approximation subspace tracking (PAST) \cite{Yang_IEEESP_1995}. GROUSE \cite{Balzano_arXiv_2010_s} proposes an incremental gradient descent algorithm performed on the Grassmannian $\mathcal{G}(d,n)$, the space of all $d$-dimensional subspace of $\mathbb{R}^n$ \cite{Edelman98a,Absil_OptAlgMatManifold_2008}. The algorithm minimizes on $\ell 2$-norm cost. GRASTA \cite{He_CVPR_2012} enhances robustness against outliers by exploiting $\ell 1$-norm cost function. PETRELS \cite{Chi_IEEETransSP_2013} calculates the underlying subspace via a discounted recursive process for each row of the subspace matrix in parallel. 

As for the {\it tensor-based} online algorithm, which is our main focus in this paper, Nion and Sidiropoulos propose an adaptive algorithm to obtain the CP decompositions \cite{Nion_IEEETransSP_2009}. Yu et al. also propose an accelerated online tensor learning algorithm (ALTO) based on Tucker decomposition \cite{Yu_ICML_2015}. However, they do not deal with a missing data presence. Mardani et al. propose an online imputation algorithm based on the CP decomposition under the presence of missing data \cite{Mardani_IEEETransSP_2015}. This considers the stochastic gradient descent (SGD) for large-scale data. This work bears resemblance to the contribution of the present paper. However, considering situations in which the subspace changes dramatically and the processing speed is sufficiently faster than data acquiring speed, a faster convergence rate algorithm per iteration is crucially important to track this change. Because it is well-known that SGD shows a slow convergence rate as the experiments described later in the paper, it is not suitable for this situation. Recently, Kasai and Mishra also proposed a novel Riemannian manifold preconditioning approach for the tensor completion problem with rank constraint \cite{kasai_ICML_2016}. 
The specific Riemannian metric allows the use of versatile framework of Riemannian optimization on quotient manifolds to develop a preconditioned  SGD algorithm. 

However, all previously described algorithms are first-order algorithms.
For that reason and because of their poor curvature approximations in ill-conditioned problems, their convergence rate can be slow. One promising approach is second-order algorithms such as stochastic quasi-Newton (QN) methods using Hessian evaluations. 
Numerous reports of the literature describe studies of stochastic versions of deterministic quasi-Newton methods  \cite{Dennis_MC_1974,Powell_NP_1976,Byrd_SJNA_1987,Nocedal_NumericalOptBook_2006} with higher scalability in the number of variables for large-scale data.
AdaGrad, which estimates the diagonal of the squared root of the covariance matrix of the gradients, was proposed \cite{Duchi_JMLR_2011}. SGD-QN exploits a {\it diagonal rescaling matrix} based on the {\it secant condition} with quasi-Newton method \cite{Bordes_JMLR_2009}. A direct extension of the deterministic BFGS in terms of using stochastic gradients and Hessian approximations is known as online BFGS \cite{Schraudolph_AISTATS_2007}. Its variants include \cite{Mokhtari_IEEETranSigPro_2014,Schraudolph_AISTATS_2007,Mokhtari_JMLR_2015,Wang_SIOPT_2017,Byrd_SIOPT_2016}. 
Overall, they achieve a higher convergence rate by exploiting curvature information of the objective function. Nevertheless, it is unclear whether those approaches operate under the subspace tracking application of interest described in this paper.

\section{Proposed online low-rank tensor subspace tracking: OLSTEC}

This section defines the problem formulation and provides the details of the proposed algorithm: OLSTEC. 

\subsection{Problem formulation}

This paper addresses the problem of the low-rank tensor completion in an online manner when the rank is a priori known or estimated. Without loss of generality, we particularly examine the third order tensor of which one order increases over time. In other words, we address $\mathbfcal{Y} \in \mathbb{R}^{L \times W \times T}$ of which third order increases infinitely. Assuming that $\mathbfcal{Y}_{i_1, i_2, i_3}$ are only known for some indices $(i_1, i_2, i_3) \in \Omega$, where $\Omega$ is a subset of the complete set of indices $(i_1, i_2, i_3)$, a general \emph{batch-based} fixed-rank tensor completion problem is formulated as 
\begin{equation}
\begin{array}{lll}
\label{Eq:CostFunction}
\displaystyle{\min_{\mathbfcal{X} \in
\mathbb{R}^{L \times W \times T}} }&   
\displaystyle{\frac{1}{2}
\| \mathcal{P}_{\Omega}(\mathbfcal{X}) - 
\mathcal{P}_{\Omega}(\mathbfcal{Y}) \|^2_F} \\
{\rm subject\ to}& {\rm rank}(\mathbfcal{X}) = R,
\end{array}
\end{equation}
where the operator $\mathcal{P}_{\Omega}(\mathbfcal{X})_{i_1, i_2, i_3} = \mathbfcal{X}_{i_1, i_2, i_3}$ if $(i_1, i_2, i_3) \in \Omega$ and $\mathcal{P}_{\Omega}(\mathbfcal{X})_{i_1, i_2, i_3}  = 0$ otherwise and (with a slight abuse of notation) $\|\cdot \|_F$ is the Frobenius norm. 
${\rm rank}(\mathbfcal{X})$ is the rank of $\mathbfcal{X}$  (see \cite{Kolda_SIAMReview_2009} for a detailed discussion of tensor rank). $R \ll \{L,W,T\}$ enforces a low-rank structure. 

Then, the problem (\ref{Eq:CostFunction}) is reformulated with $\ell_2$-norm regularizers as \cite{Mardani_IEEETransSP_2015}
\begin{equation}  
\label{Eq:Batch_Problem_Definition}
\min_{\scriptsize{\mat{A},\mat{B},\mat{C}}} 
\frac{1}{2} \| \mathcal{P}_{\Omega}(\mathbfcal{Y}) -  \mathcal{P}_{\Omega}(\mathbfcal{X}) \|_F^2 + 
\mu (\| \mat{A}\|_F^2 + \| \mat{B}\|_F^2 + \| \mat{C}\|_F^2) \nonumber
\end{equation}
\begin{equation}  
 {\rm subject\ to} \ \ \ \mat{X}_\tau =  \mat{A} {\rm diag}(\vec{b}^{\tau}) \mat{C}^T   \ \ \ \ \ \ 
{\rm \ for} \ \tau = 1, \ldots, T.
\end{equation}  
where $\mu>0$ is a regularizer parameter. This regularizer suppresses the instability of RLS described in Section \ref{Sec:RLS}. Consequently, considering the situation where the partially observed tensor slice ${\bf \Omega}_\tau \circledast \mat{Y}_\tau$ is acquired sequentially over time, we estimate $\{ \mat{A}, \mat{B}, \mat{C}\}$ by minimizing the exponentially weighted least squares;
        \begin{eqnarray}
        	\label{Eq:Final_Problem_Definition}
        \min_{\scriptsize{\mat{A},\mat{B},\mat{C}}}  \frac{1}{2} \sum_{\tau=1}^t \lambda^{t-\tau} 
        \biggl[ {\| {\bf \Omega}_\tau  \circledast \bigl[ \mat{Y}_\tau -  \mat{A} {\rm diag}(\vec{b}^\tau) \mat{C}^T 
        \bigr] \|_F^2}  +\ {\bar{\mu}(\| \mat{A}\|_F^2 + \| \mat{C}\|_F^2) + \mu_r[\tau] \| \vec{b}^\tau \|_2^2 } 
\biggr].
        \end{eqnarray}
Therein, $\mu_r[t]$ is the regularizer parameter for \vec{b}, $\bar{\mu}=\mu_r[\tau]/\sum_{\tau=1}^t \lambda^{t-\tau}$, and $0 < \lambda \leq 1$ is the so-called forgetting parameter. $\lambda=1$ case is equivalent to the batch-based problem (\ref{Eq:Batch_Problem_Definition}).

\subsection{Algorithm of OLSTEC}
\label{Sec:AlgorithmOfOLSTEC}
The unknown variables in (\ref{Eq:Final_Problem_Definition})  are $\mat{A}, \mat{C}$, and $\vec{b}$. Also, \mat{A} and \mat{C} are a non-convex set. Therefore, this function is non-convex. The proposed OLSTEC algorithm, as summarized by Algorithm 1, alternates between least-squares estimation of $\vec{b}[t]$ for fixed $\mat{A}[t\!-\!1]$ and $\mat{C}[t\!-\!1]$, and a second-order stochastic gradient step using the RLS algorithm on $\mat{A}[t\!-\!1]$ and $\mat{C}[t\!-\!1]$ for fixed $\vec{b}[t]$. 
It is noteworthy, as described in Section \ref{sec:notations}, that $\mat{W}[t]$ with the square bracket represents the calculated $\mat{W}$ after performing $t$-times updates.

\subsubsection{Calculation of $\vec{b}[t]$}

%
%

The estimate $\vec{b}[t]$ is obtained in a closed form by minimizing the residual by fixing $\{ \mat{A}[t-1], \mat{C}[t-1] \}$ derived at time $t-1$. Hereinafter, in this subsection, we use the notation of $\mu_r$ for $\mu_r[t]$.
\begin{eqnarray}
\label{AppenEq:Problem_Definition_b_timebase}
\vec{b}[t] & = & \defargmin_{\scriptsize \vec{b} \in \mathbb{R}^R} \frac{1}{2} 
\Biggl[ \|  {\bf \Omega}_t \circledast 
[ \mat{Y}_t -  \mat{A}[t\!-\!1] {\rm diag}(\vec{b}) (\mat{C}[t\!-\!1])^T ] \|_F^2 + \mu_r \| \vec{b} \|_2^2 \biggr] 
 \nonumber \\
  & = & 
  \defargmin_{\scriptsize \vec{b} \in \mathbb{R}^R} \frac{1}{2} 
\biggl[  \sum_{(l,w) \in {\bf \Omega}_t} \left( [\mat{Y}_t]_{l,w} - (\vec{a}^l[t\!-\!1] \circledast \vec{c}_w[t\!-\!1])^T \vec{b} \right)^2  + \mu_r \| \vec{b} \|_2^2 \biggr] \nonumber \\
& = & \defargmin_{\scriptsize \vec{b} \in \mathbb{R}^R} 
\frac{1}{2} 
\Biggl[ \sum_{l=1}^L \sum_{w=1}^W 
\left(
[{\bf \Omega}_t]_{l,w}  
\left(
[\mat{Y}_t]_{l,w} -
 (\vec{g}_{l,w}[t] )^T  \vec{b}
  \right) \right)^2 
+ \mu_r \| \vec{b} \|_2^2
\Biggr].
\end{eqnarray}
Therein, $\vec{g}_{l,w}[t]= \vec{a}^{l}[t\!-\!1] \circledast \vec{c}^{w}[t\!-\!1] \in \mathbb{R}^R$. Defining $F[t]$ as the inner objective to be minimized, we obtain the following. 
\begin{eqnarray*} 
%
\frac{\partial F[t](\vec{b})}{\partial \vec{b}} & = & 
- \sum_{l=1}^L \sum_{w=1}^W 
[{\bf \Omega}_t]_{l,w}  
\left(
\mat{Y}[t]_{l,w} -
  (\vec{g}_{l,w}[t] )^T  \vec{b}
\right)  \vec{g}_{l,w}[t] 
+ \mu_r \vec{b}.
\end{eqnarray*} 
Because $\vec{b}$ satisfies $\partial F[t]/\partial \vec{b} = 0$, we obtain $\vec{b}$ as shown below.
\begin{eqnarray*} 
\mu_r \vec{b} +  
\sum_{l=1}^L  \sum_{w=1}^W [{\bf \Omega}_t]_{l,w} 
  (\vec{g}_{l,w}[t] )^T  \vec{b} \ \vec{g}_{l,w}[t] 
=  
\sum_{l=1}^L \sum_{w=1}^W 
[{\bf \Omega}_t]_{l,w}  \mat{Y}[t]_{l,w}
\vec{g}_{l,w}[t].
\end{eqnarray*} 

Finally, we obtain $\vec{b}[t]$ as
\begin{eqnarray}
\label{Eq:Solution_b_timebase}
\vec{b}[t]  =  
\biggl[ 
\mu_r \mat{I}_R +  
\sum_{l=1}^L  \sum_{w=1}^W [{\bf \Omega}_t]_{l,w} 
 \vec{g}_{l,w}[t]  (\vec{g}_{l,w}[t] )^T\biggr]^{-1}
\biggl[
\sum_{l=1}^L \sum_{w=1}^W 
[{\bf \Omega}_t]_{l,w}  \mat{Y}[t]_{l,w} 
\vec{g}_{l,w}[t]
\biggr].
\end{eqnarray}   

\begin{algorithm}[t]
\caption{OLSTEC algorithm}
\label{alg:algorithm}
\begin{algorithmic}[1]
\REQUIRE{ $\{ \mat{Y}_t$ and ${\bf {\Omega}}_t \}^{\infty}_{t=1}$, $\lambda$, $\mu$}
\STATE{Initialize \{$\mat{A}[0]$, $\vec{b}[0]$, $\mat{C}[0]$\}, $\mat{Y}[0]=\mat{0}$, $(\mat{RA}_l[0])^{-1}=(\mat{RC}_w[0])^{-1}=\gamma \mat{I}_{R}, \gamma > 0$.}
\FOR{$t=1,2, \cdots$} 
	\STATE{Calculate $\vec{b}[t]$ \hfill Equation (\ref{Eq:Solution_b_timebase})}
	\STATE{$\mat{X}_t = \mat{A}[t\!-\!1] {\rm diag}(\vec{b}_t) (\mat{C}[t\!-\!1])^{T}$}
	\FOR{$l=1,2, \cdots, L$} 
	\STATE{Calculate $\mat{RA}_l[t]$ \hfill Equation (\ref{Eq:Update_RA})}
	\STATE{Calculate $\vec{a}^l[t]$ \hfill Equation (\ref{Eq:al_final})}
	\ENDFOR
	\FOR{$w=1,2, \cdots, W$} 
	\STATE{Calculate $\mat{RC}_l[t]$ \hfill Equation (\ref{Eq:Update_RC})}
	\STATE{Calculate $\vec{c}^w[t]$ \hfill Equation (\ref{Eq:cw_final})}
	\ENDFOR	
\ENDFOR
\RETURN $\mat{X}_t = \mat{A}[t] {\rm diag}(\vec{b}[t]) (\mat{C}[t])^{T}$
\end{algorithmic}
\end{algorithm}

\subsubsection{Calculation of $\mat{A}[t]$ and $\mat{C}[t]$ based on RLS}

The calculation of $\mat{C}[t]$ uses $\mat{A}[t\!-\!1]$. The calculation of $\mat{A}[t]$ uses $\mat{C}[t\!-\!1]$. This paper addresses a second-order stochastic gradient based on the RLS algorithm with forgetting parameters, which has been used widely in tracking of time varying parameters in many fields. Its computation is efficient because we update the estimates recursively every time new data become available.

As for $\mat{A}[t]$, the problem (\ref{Eq:Final_Problem_Definition}) is reformulated as
\begin{eqnarray}
\label{Eq:Problem_Definition_A_timebase}
\min_{\scriptsize \mat{A} \in \mathbb{R}^{L \times R}} \frac{1}{2} \sum_{\tau}^t \lambda^{t-\tau} 
\biggl[ \| {\bf \Omega}_\tau \circledast \bigl[ \mat{Y}_\tau - \mat{A} {\rm diag}(\vec{b}[\tau]) (\mat{C}[\tau\!-\!1])^T 
\bigr]  \|_F^2 \biggr] + \frac{\mu_r[t]}{2} \| \mat{A}\|_F^2.
\end{eqnarray}

The objective function in (\ref{Eq:Problem_Definition_A_timebase}) decomposes into a parallel set of smaller problems, one for each row of 
$\mat{A}$, as
\begin{eqnarray}
	\label{Eq:problem_def_am}
	\vec{a}^l[t] &  = & \defargmin_{\vec{a}^l \in \mathbb{R}^{R}} \frac{1}{2}  
	\sum_{\tau=1}^t 
	\Biggl[ \lambda^{t-\tau}
	\sum_{w=1}^W
	  [{\bf \Omega}_\tau]_{l,w}\left(
	  [\mat{Y}_\tau]_{l,w} -(\vec{a}^l)^T  
	  {\rm diag} (\vec{b}[\tau])
	\vec{c}^w[\tau\!-\!1] \right) ^2 
	\Biggr] \nonumber \\
 &&+ \frac{\mu_r[t]}{2} \| \vec{a}^l\|_2^2. \ \ \ \ \ \ \ \ 
\end{eqnarray}

Here, denoting ${\rm diag}(\vec{b}[\tau]) \vec{c}^w[\tau\!-\!1]$ as $\vec{\alpha}_w[\tau] \in \mathbb{R}^R$, the derivative of the inner objective function in (\ref{Eq:problem_def_am}) with regard to $\vec{a}^l \in \mathbb{R}^{R}$ is calculated as shown below.
\begin{eqnarray}
	\frac{\partial{F(\vec{a}^l)}}{\partial{(\vec{a}^l})} & = & 
	\sum_{\tau=1}^t 
	\biggl[
	-
	\sum_{w=1}^W  
	\lambda^{t-\tau} [{\bf \Omega}_\tau]_{l,w}
	\left([\mat{Y}_\tau]_{l,w} 
	- (\vec{a}^l)^T \vec{\alpha}_w[\tau] \right)  \vec{\alpha}_{w} [\tau] \biggr]
	+ \mu_r[t]\vec{a}^l.\nonumber
\end{eqnarray}

Then, by setting this derivative equal to zero, we get $\vec{a}^l[t]$ as
\begin{equation}
\begin{split}
	%
	&\left( \sum_{\tau=1}^t 
	\biggl[
	\sum_{w=1}^W  \lambda^{t-\tau} [{\bf \Omega}_\tau]_{l,w} 
	\vec{\alpha}_w[\tau] (\vec{\alpha}_w[\tau])^T 
	\biggr]+ \mu_r[t]\mat{I}_R \right)
	\vec{a}^l[t]\nonumber
	\\
	&\hspace*{2cm}= \sum_{\tau=1}^t 
	\sum_{w=1}^W  \lambda^{t-\tau}  [{\bf \Omega}_\tau]_{l,w} [\mat{Y}_\tau]_{l,w} 
	 \vec{\alpha}_w[\tau].
\end{split}	 
\end{equation}
Therein, $(\vec{a}^l)^T \vec{\alpha}_w[\tau] \vec{\alpha}_w[\tau]= (\vec{\alpha}_w[\tau])^T \vec{a}^l \vec{\alpha}_w[\tau]
= \vec{\alpha}_w[\tau](\vec{\alpha}_w[\tau])^T \vec{a}^l$ is used. Finally, we obtain the following as
\begin{eqnarray}	 
	\label{Eq:RA_a_s}
	\mat{RA}_l[t]  \vec{a}^l[t]  & = &  \vec{s}_l[t],\nonumber
\end{eqnarray}	
where $\mat{RA}_l[t] \in \mathbb{R}^{R\times R}$ and $\vec{s}_l[t] \in \mathbb{R}^{R}$ are defined as shown below.
\begin{eqnarray}	
 \mat{RA}_l[t] &  =  &   \sum_{\tau=1}^t 
 	\biggl[
 	\sum_{w=1}^W  \lambda^{t-\tau} [{\bf \Omega}_\tau]_{l,w} 
	\vec{\alpha}_w[\tau] \vec{\alpha}_w[\tau]^T
	\biggr] + \mu_r[t] \mat{I}_R	\nonumber \\
	\label{Eq:s_l}
	\vec{s}_l[t] &  =   &  \sum_{\tau=1}^t 
	\biggl[
	\sum_{w=1}^W  \lambda^{t-\tau}  [{\bf \Omega}_\tau]_{l,w}  [\mat{Y}_\tau]_{l,w} 
	  \vec{\alpha}_w[\tau]
	 \biggr].\nonumber
\end{eqnarray}	
Here, $\mat{RA}_l[t]$ is transformed as
\begin{eqnarray}
	\label{Eq:Update_RA}
	\mat{RA}_l[t]  & = & 
	\sum_{\tau=1}^{t-1}
 	\left(
 	\sum_{w=1}^W  \lambda^{t-\tau} [{\bf \Omega}_\tau]_{l,w} 
	\vec{\alpha}_w[\tau] (\vec{\alpha}_w[\tau])^T
	\right)
+ \sum_{w=1}^W [{\bf \Omega}_t]_{l,w}  \vec{\alpha}_w[t]  (\vec{\alpha}_w[t])^T
	+ \mu_r[t]  \mat{I}_R
 	\nonumber \\
	&=& \lambda 
	\sum_{\tau=1}^{t-1}
 	\left(
 	\sum_{w=1}^W  \lambda^{t-1-\tau} [{\bf \Omega}_\tau]_{l,w} 
	\vec{\alpha}_w[\tau] (\vec{\alpha}_w[\tau])^T
	\right)	
 +  \sum_{w=1}^W [{\bf \Omega}_t]_{l,w}  \vec{\alpha}_w[t]  (\vec{\alpha}_w[t])^T
	+ \mu_r[t]  \mat{I}_R \nonumber \\
	&=& \lambda \biggl[
	\sum_{\tau=1}^{t-1}
 	\left(
 	\sum_{w=1}^W  \lambda^{t-1-\tau} [{\bf \Omega}_\tau]_{l,w} 
	\vec{\alpha}_w[\tau] (\vec{\alpha}_w[\tau])^T \right)
	+ \mu_r[\tau]  \mat{I}_R  \biggr]
	\nonumber \\
	&&+  \sum_{w=1}^W [{\bf \Omega}_t]_{l,w}  \vec{\alpha}_w[t]  (\vec{\alpha}_w[t])^T
 	+ \mu_r[t]  \mat{I}_R - \lambda \mu_r[t\!-\!1]  \mat{I}_R \nonumber \\
 	&=&  \lambda \mat{RA}_l[t\!-\!1]
	+  \sum_{w=1}^W [{\bf \Omega}_t]_{l,w}  \vec{\alpha}_w[t]  (\vec{\alpha}_w[t])^T
	+ (\mu_r[t]  - \lambda \mu_r[t\!-\!1] ) \mat{I}_R.
\end{eqnarray}
%
%

Similarly, $\vec{s}_l[t]$ is also transformed as 
\begin{eqnarray}
	\label{Eq:Update_s}
	\vec{s}_l[t] 
	& = & 
	\sum_{w=1}^W \biggl[  
	\lambda^{t-1}  [{\bf \Omega}_1]_{l,w}  [\mat{Y}_{1}]_{l,w} \vec{\alpha}_{w}[1] 
	\cdots +  \lambda^{1} [{\bf \Omega}_{t\!-\!1}]_{l,w}   [\mat{Y}_{t\!-\!1}]_{l,w}   \vec{\alpha}_{w}[t\!-\!1] \nonumber \\
	&&+ \lambda^{0} [{\bf \Omega}_t]_{l,w}  [\mat{Y}_t]_{l,w} \vec{\alpha}_w[t]  \biggr]\nonumber \\
	& = & \lambda \vec{s}_l[t\!-\!1] + \sum_{w=1}^W [{\bf \Omega}_t]_{l,w} [\mat{Y}_t]_{l,w} \vec{\alpha}_w[t]. \nonumber
\end{eqnarray}

%
From $\mat{RA}_l[t] \vec{a}^l[t]  = \vec{s}_l[t] $, we modify the above as shown below.
%
\begin{eqnarray}
	\label{Eq:RA_a_update}
	\mat{RA}_l[t] \vec{a}^l[t]  
	& = &  \lambda \vec{s}_l[t\!-\!1] + \sum_{w=1}^W [{\bf \Omega}_t]_{l,w}  [\mat{Y}_t]_{l,w}  \vec{\alpha}_w[t] \nonumber \\
	& = &  \lambda \mat{RA}_l[t\!-\!1]  \vec{a}^l[t\!-\!1] 
	 +\sum_{w=1}^W [{\bf \Omega}_t]_{l,w}   [\mat{Y}_t]_{l,w}   \vec{\alpha}_w[t] 
	\nonumber \\
	& = &  \left(\mat{RA}_l[t] - \sum_{w=1}^W [{\bf \Omega}_t]_{l,w}  \vec{\alpha}_w[t]  (\vec{\alpha}_w[t])^T 
	- (\mu_r[t] \!-\! \lambda \mu_r[t\!-\!1]) \mat{I}_R\right)		 \vec{a}^l[t\!-\!1]  \nonumber\\
	 &&+ 
	\sum_{w=1}^W [{\bf \Omega}_t]_{l,w}  [\mat{Y}_t]_{l,w}   \vec{\alpha}_w[t] 
	\nonumber \\
	& = &  \left(\mat{RA}_l[t] 
	 - (\mu_r[t] - \lambda \mu_r[t\!-\!1]) \mat{I}_R\right) 
	 \vec{a}^l[t\!-\!1]  \nonumber \\
	&&+ 
	\sum_{w=1}^W [{\bf \Omega}_t]_{l,w} \left([\mat{Y}_t]_{l,w} - (\vec{\alpha}_w[t])^T  \vec{a}^l[t\!-\!1] \right) \vec{\alpha}_w[t].
\end{eqnarray}


Subsequently, $\vec{a}^l[t]$ is obtained  as presented below.
\begin{eqnarray}
	\label{Eq:al_final}
	\vec{a}^l[t] &=& \vec{a}^l[t\!-\!1]  - (\mu[t] - \lambda \mu[t\!-\!1]) (\mat{RA}_l[t])^{-1}\vec{a}^l[t\!-\!1]  \nonumber  \\
	 &&+ 
	 \sum_{w=1}^W [{\bf \Omega}_{t}]_{l,w} 
	  \left([\mat{Y}_{t}]_{l,w}  -  (\vec{\alpha}_w[t])^T \vec{a}^l[t\!-\!1] \right)  
	(\mat{RA}_l[t])^{-1} \vec{\alpha}_w[t]. \hspace*{1cm}
\end{eqnarray}
As in the $\mat{A}[t]$ case, $\vec{c}^w[t] \in \mathbb{R}^{R}$ is obtainable by denoting $(\vec{a}^l[\tau])^T {\rm diag}(\vec{b}[{\tau}])$ as $\vec{\beta}^{l}[\tau] \in \mathbb{R}^{1 \times R}$ as 
\begin{eqnarray}
	\label{AppenEq:derivative_by_cw}
	\frac{\partial{G(\vec{c}^w)}}{\partial{(\vec{c}^w})} & = & 
	\sum_{\tau=1}^t 
	\biggl[
	\sum_{l=1}^L  
	\lambda^{t-\tau} [{\bf \Omega}_\tau]_{l,w}
	\left([\mat{Y}_\tau]_{l,w} - \vec{\beta}^{l}[\tau] \vec{c}^w \right) (-\vec{\beta}^{l} [\tau])^T
	\biggr] + \mu_r[t]\vec{c}^w.\hspace*{1cm}\nonumber
\end{eqnarray}

Then, we obtain $\vec{c}^w$ as presented below.
\begin{eqnarray}
	\label{}
	%
	\left(
	\sum_{\tau=1}^t 
	\sum_{l=1}^L  
	\lambda^{t-\tau} [{\bf \Omega}_\tau]_{l,w}
	\vec{\beta}^{l} [\tau] (\vec{\beta}^{l} [\tau])^T
	+ \mu_r[t] \mat{I}_R
	\right)\vec{c}^w 
	= 
	\sum_{\tau=1}^t 
	\sum_{l=1}^L
	\lambda^{t-\tau} [{\bf \Omega}_\tau]_{l,w}
	[\mat{Y}_\tau]_{l,w} 
	(\vec{\beta}^{l} [\tau])^T.	\nonumber
\end{eqnarray}

Finally, we obtain the following.
\begin{eqnarray}	 
	\mat{RC}_w[t] \vec{c}^w[t]  & = &  \vec{s}_w[t] \nonumber \\
	\vec{c}^w[t]  & =  & (\mat{RC}_w[t])^{\dagger} \vec{s}_w[t].\nonumber
\end{eqnarray}	
In that expression, $\mat{RC}_w[t] \in \mathbb{R}^{R\times R}$ and $\vec{s}_w[t] \in \mathbb{R}^{R}$ are defined as 
\begin{eqnarray}	
 \mat{RC}_w[t] & :=& 	
 	\sum_{\tau=1}^t 
	\sum_{l=1}^L  
	\lambda^{t-\tau} [{\bf \Omega}_\tau]_{l,w}
	(\vec{\beta}^{l} [\tau])^T \vec{\beta}^{l} [\tau]
	+ \mu_r[t] \mat{I}_R, \nonumber \\
	\vec{s}_w[t] & := &  	\sum_{\tau=1}^t 
	\sum_{l=1}^L
	\lambda^{t-\tau} [{\bf \Omega}_\tau]_{l,w}
	[\mat{Y}_\tau]_{l,w}(\vec{\beta}^{l} [\tau])^T.\nonumber
\end{eqnarray}	

Here, $\mat{RC}_w[t]$ is transformed as 
\begin{eqnarray}
	\label{Eq:Update_RC}
	\mat{RC}_w[t]  & = & 
 	\sum_{\tau=1}^{t-1} 
	\sum_{l=1}^L  
	\lambda^{t-\tau} [{\bf \Omega}_\tau]_{l,w}
	(\vec{\beta}^{l} [\tau])^T\vec{\beta}^{l} [\tau] 
	+\sum_{l=1}^L  
	[{\bf \Omega}_t]_{l,w}\vec{\beta}^{l} [t] (\vec{\beta}^{l} [t])^T
	+ \mu_r[t] \mat{I}_R \nonumber \\
	& = & 
 	\lambda \mat{RC}_w[t-1] 
	+\sum_{l=1}^L  
	[{\bf \Omega}_t]_{l,w} (\vec{\beta}^{l} [t])^T \vec{\beta}^{l} [t] 	
	+ (\mu_r[t] - \lambda \mu_r[t\!-\!1]) \mat{I}_R. \nonumber
\end{eqnarray}

Similarly, $\vec{s}_w[t]$ is also transformed as shown below.
\begin{eqnarray}
	\label{AppenEq:Update_sw}
	\vec{s}_w[t] 
	& = & \lambda \vec{s}_w[t\!-\!1] + \sum_{l=1}^L [{\bf \Omega}_t]_{l,w}[\mat{Y}_t]_{l,w}(\vec{\beta}^{l} [t])^T. \nonumber
\end{eqnarray}

From $\mat{RC}_w[t] \vec{c}^w[t]  = \vec{s}_w[t] $, we can modify the following as
\begin{eqnarray}
	\label{AppenEq:RC_a}
	\mat{RC}_w[t] \vec{c}^w[t]  
	& = & 
	\left(
	\mat{RC}_w[t] 	
	+ (\mu_r[t]  - \lambda \mu_r[t-1])\mat{I}_R
	\right) \vec{c}^w[t\!-\!1] \nonumber \\
	&&+ \sum_{l=1}^L[{\bf \Omega}_t]_{l,w}
	\left(
	[\mat{Y}_t]_{l,w} - \vec{\beta}^{l}[t] \vec{c}^w[t\!-\!1]\right)(\vec{\beta}^{l}[t])^T. 
\end{eqnarray}
Finally, $\vec{c}^w[t] \in \mathbb{R}^{R}$ is obtained as
\begin{eqnarray}
	\label{Eq:cw_final}
	\vec{c}^w[t] 
	 & = & \vec{c}^w[t\!-\!1] - (\mu_r[t]  - \lambda \mu_r[t-1]) (\mat{RC}_w[t])^{-1} \vec{c}^w[t\!-\!1] \nonumber \\
	 &&+ 
	 \sum_{l=1}^L[{\bf \Omega}_t]_{l,w}
	\left(
	[\mat{Y}_t]_{l,w} - \vec{\beta}^{l}[t] \vec{c}^w[t\!-\!1]  \right) (\mat{RC}_w[t])^{-1} (\vec{\beta}^{l}[t])^T. \hspace*{1cm}
\end{eqnarray}


\section{Extensions of the OLSTEC algorithm}	

Two extensions are explored in this section. 

\subsection{Truncated window setting}
\label{Sec:TruncatedWindowSetting}

Given a truncated window of length $V \leq T$, we deal with a truncated tensor $\mathbfcal{Y}$ of size 
$L \times W \times V$. Without ignoring the regularizer for simplicity, the sub problem (\ref{Eq:problem_def_am}) is 
reformulated as shown below.
\begin{eqnarray}
	\label{Eq:problem_def_al_TW}
	\vec{a}^l[t]  = \defargmin_{\vec{a}^l}  \frac{1}{2}  \sum_{\tau=t-V+1}^t \biggl[ 
	\lambda^{t-\tau}  \sum_{w=1}^W [{\bf \Omega}_{\tau}]_{l,w} 
	( [{\bf Y}_{\tau}]_{l,w}  
	- (\vec{a}^l)^T {\rm diag} (\vec{b}[\tau])\vec{c}^w[\tau-1] )^2 \biggr].\nonumber
\end{eqnarray}

The update of $\mat{RA}_l[t]$ in (\ref{Eq:Update_RA}) and $\vec{s}_l[t]$ in (\ref{Eq:Update_s}) are modified, respectively, as presented below.
\begin{eqnarray}
	\label{Eq:TW_RA_s}
	\mat{RA}_l[t]  & = &  \displaystyle{\lambda \mat{RA}_l^{[t-1]} + \sum_{w=1}^W [{\bf \Omega}_{t}]_{l,w}  \vec{\alpha}_{w}[t]  \vec{\alpha}_{w}[t]^T} \nonumber\\
	&&\displaystyle{- \sum_{w=1}^W \lambda^{V} [{\bf \Omega}_{t-V}]_{l,w}  \vec{\alpha}_w[t-V]  \vec{\alpha}_w[t-V]^T}\nonumber\\
	\vec{s}_l[t] & = &  \displaystyle{\lambda \vec{s}_l[t-1] + \sum_{w=1}^W [{\bf \Omega}_{t}]_{l,w}   [{\bf Y}_{t}]_{l,w}   \vec{\alpha}_{w}[t]} \nonumber\\
	&&\displaystyle{-  \sum_{w=1}^W \lambda^{V} [{\bf \Omega}_{t-V}]_{l,w} [{\bf Y}_{t-V}]_{l,w}  \vec{\alpha}_w[t-V]}.\nonumber
\end{eqnarray}

Therefore, the statement below corresponds to (\ref{Eq:RA_a_update}).
\begin{eqnarray}
	\mat{RA}_l[t] \cdot  \vec{a}^l[t] 
	& = & \lambda \vec{s}_l[t-1] + \sum_{w=1}^W [{\bf \Omega}_{t}]_{l,w}   [{\bf Y}_{t}]_{l,w}   \vec{\alpha}_{w}[t] 
	-  \sum_{w=1}^W \lambda^{V} [{\bf \Omega}_{t-V}]_{l,w}  [{\bf Y}_{t-V}]_{l,w}   \vec{\alpha}_w[t-V]
	\nonumber \\
	%
	& = & \mat{RA}_l[t] \vec{a}^l[t-1]  + \sum_{w=1}^W [{\bf \Omega}_{t}]_{l,w}  \biggl[[{\bf Y}_{t}]_{l,w}  
	-  (\vec{\alpha}_{w}[t])^T \vec{a}^l[t-1]  \biggr]\vec{\alpha}_{w}[t] \nonumber \\
	& & - \sum_{w=1}^W \lambda^{V} [{\bf \Omega}_{t-V}]_{l,w} \biggl[ [{\bf Y}_{t-V}]_{l,w}  -  (\vec{\alpha}_w[t-V])^T  \vec{a}^l[t-1] \biggr] \vec{\alpha}_w[t-V]. \nonumber
\end{eqnarray}
Finally, $\vec{a}^l[t]$ in (\ref{Eq:al_final}) is replaced with the following.
\begin{eqnarray}
	\label{Eq:am_final_truncated_window}
	\vec{a}^l[t] & = & \vec{a}^l[t-1]  + \sum_{w=1}^W [{\bf \Omega}_{t}]_{l,w}  \biggl[[{\bf Y}_{t}]_{l,w}  
	-  (\vec{\alpha}_{w}[t])^T \vec{a}^l[t-1]  \biggr]  (\mat{RA}_l[t])^{-1} \vec{\alpha}_{w}[t]\nonumber \\
	& & - \sum_{w=1}^W \lambda^{V} [{\bf \Omega}_{t-V}]_{l,w} \biggl[ [{\bf Y}_{t-V}]_{l,w}  -  (\vec{\alpha}_w[t-V])^T  \vec{a}^l[t-1]  \biggr] 
	 (\mat{RA}_l[t])^{-1} \vec{\alpha}_w[t-V].\nonumber
\end{eqnarray}
Similarly, instead of (\ref{Eq:cw_final}), $\vec{c}^w[t]$ is obtainable as shown below.
\begin{eqnarray}
	\label{Eq:cn_final_truncated_window}
	\vec{c}^w[t] & = & \vec{c}^w[t-1]  + \sum_{l=1}^L [{\bf \Omega}_{t}]_{l,w}  \biggl[[{\bf Y}_{t}]_{l,w}  
	-   \vec{\beta}^l[t] \vec{c}^w[t-1]  \biggr]  (\mat{RC}_w[t])^{-1} (\vec{\beta}^l[t])^T \nonumber \\
	& & - \sum_{l=1}^L \lambda^{V} [{\bf \Omega}_{t-V}]_{l,w} \biggl[ [{\bf Y}_{t-V}]_{l,w}  - \vec{\beta}^l[t-V] \vec{c}^w[t-1] \biggr] 
	 (\mat{RC}_w[t])^{-1} 
	 (\vec{\beta}^l[t-V])^T.\nonumber
\end{eqnarray}

\subsection{Simplified OLSTEC}
\label{Sec:SimplifiedOLSTEC}

The calculation costs of the inversions of $\mat{RA}[t]$ and $\mat{RC}[t]$ are the most expensive parts in (\ref{Eq:al_final}) and (\ref{Eq:cw_final}). Therefore, we execute an diagonal approximation of $\mat{RA}[t]$ and $\mat{RC}[t]$ to reduce the calculation costs, which ignores the off-diagonal part of them. More specifically, we calculate $\vec{a}^l[t]$ instead of (\ref{Eq:al_final}) as presented below.

\begin{eqnarray}
	\label{Eq:al_final_simplified}
	\vec{a}^l[t] &=& \vec{a}^l[t\!-\!1]  - (\mu[t] - \lambda \mu[t\!-\!1]) (\mat{DA}_l[t])^{-1}\vec{a}^l[t\!-\!1]  \nonumber  \\
	 &&+ 
	 \sum_{w=1}^W [{\bf \Omega}_{t}]_{l,w} 
	  \left([\mat{Y}_{t}]_{l,w}  -  (\vec{\alpha}_w[t])^T \vec{a}^l[t\!-\!1] \right)  
	(\mat{DA}_l[t])^{-1} \vec{\alpha}_w[t].\nonumber
\end{eqnarray}
Therein, $\mat{DA}_l[t]$ is defined by reformulating (\ref{Eq:Update_RA}) as
\begin{eqnarray}
	\label{Eq:DA}
	\mat{DA}_l[t] =  \lambda \mat{DA}_l[t\!-\!1]
	+  {\rm diag}\left(\sum_{w=1}^W [{\bf \Omega}_t]_{l,w}  \vec{\alpha}_w[t]  (\vec{\alpha}_w[t])^T \right)
	+ (\mu_r[t]  - \lambda \mu_r[t\!-\!1] ) \mat{I}_R.
\end{eqnarray}
The calculation of $(\mat{DA}_l[t])^{-1}$ is very light because $\mat{DA}_l[t]$ is a diagonal matrix. Similarly, the calculation of $\vec{c}^w[t]$ can be simplified. 
\section{Computational complexity and memory consumption analysis}

With respect to computational complexity per iteration, 
OLSTEC requires $\mathcal{O}(| {\bf {\Omega}}_t | R^2+(L+W) R^3)$ because of $\mathcal{O}(| {\bf {\Omega}}_t | R^2)$ for $\vec{b}[t]$ in (\ref{Eq:Solution_b_timebase}) and $\mathcal{O}((L+W) R^3)$ for the inversion of $\mat{RA}_l$ and $\mat{RC}_w$ in (\ref{Eq:al_final}) and (\ref{Eq:cw_final}), respectively. 
As for memory consumption, $\mathcal{O}((L+W)R^2)$ is required, respectively, for $\mat{RA}_l$ and $\mat{RC}_w$. 
The simplified OLSTEC requires $\mathcal{O}(| {\bf {\Omega}}_t | R^2+(L+W)R)$, where the second term is achieved by the diagonal approximation of $\mat{RA}_l$ and $\mat{RC}_w$ such as (\ref{Eq:DA}). Similarly, the memory consumption is reduced to $\mathcal{O}((L+W)R)$.

\section{Numerical Evaluations}

We present numerical comparisons of the OLSTEC algorithm with state-of-the-art algorithms on synthetic and real-world datasets. All the following experiments are done on a PC with 2.6 GHz Intel Core i7 CPU and 16 GB RAM. We implement our proposed algorithm in Matlab\footnote{The Matlab code of OLSTEC is available in \url{https://github.com/hiroyuki-kasai/OLSTEC}}, and use Matlab codes of the comparing algorithm provided by the respective authors except for the algorithm, designated as  ``TeCPSGD" algorithm in this paper, proposed in \cite{Mardani_IEEETransSP_2015}.  Finally, as for evaluation metrics, we use {\it the normalized residual error} at each iteration $t$ defined as
\begin{eqnarray}
\text{Normalized residual error} & : & \frac{\| \mat{X}_t - \mat{Y}_t\|^2_F}{\| \mat{Y}_t \|^2_F},\nonumber
\end{eqnarray}
and {\it the running averaging error} defined as
\begin{eqnarray}
\text{Running averaging error} & : &\frac{1}{T} \sum_{\tau=1}^T \frac{\| \mat{X}_\tau - \mat{Y}_\tau\|^2_F}{\| \mat{Y}_\tau \|^2_F}.\nonumber
\end{eqnarray}

\subsection{Synthetic datasets}

We first evaluate the performance of our proposed algorithm using synthetic datasets compared with TeCPSGD. The two versions of the proposed algorithm, i.e., the original OLSTEC in Section \ref{Sec:AlgorithmOfOLSTEC} and the simplified OLSTEC in Section\ref{Sec:SimplifiedOLSTEC}, are also compared. This experiment specifically considers a time-varying subspace to demonstrate the capability of the proposed algorithm to handle it. To this end, we generate the synthetic low rank-$R$ tensors $\mathbfcal{Y} \in \mathbb{R}^{L\times W \times T}$ with the factor matrices $\mat{A}$ and $\mat{C}$ by updating them in the direction of the third order direction as $\mat{A}_{t+1}=\mat{A}_{t} \mat{Q}(t, \alpha)$ and $\mat{C}_{t+1}=\mat{C}_{t} \mat{Q}(t, \alpha)$.  $\mat{Q}(t, \alpha) \in \mathbb{R}^{R \times R}$ is the rotation matrix given as 
\begin{eqnarray}
\mat{Q}(t, \alpha)&=&\left(
\begin{array}{c|cc|c}
  \mat{I}_{p-1} & 0 & 0& 0\\ 
  \hline
  0 & \cos(\alpha)& -\sin(\alpha)& 0\\ 
  0 & \sin(\alpha) & \cos(\alpha)& 0\\   
  \hline
  0 & 0 & 0& \mat{I}_{R-p-1}\\ 
\end{array}
\right),\nonumber
\end{eqnarray}
where $p(t) = (t + R-2)\%(R-1) + 1$, and $\alpha$ is the rotation angle. The higher values of $\alpha$ mimic more dynamically changing subspace. As the definition of $\mat{Q}(t, \alpha)$ clearly represents, the direction of the rotation is changed every iteration $t$. $\mat{A}_1$ and $\mat{C}_1$ are generated with i.i.d standard Gaussian $\mathcal{N}(0, 1)$ entries. Also, Gaussian noise with i.i.d $\mathcal{N}(0, \epsilon^2)$ entries are added. We set $T=500$, $R=5$. The noise level is $\epsilon=10^{-3}$. $\mu_r[t]=10^{-3}$ and $\lambda=0.5$ are set in the proposed algorithm. $\lambda$, $\mu$ and the stepsize are set, respectively, to $0.001$, $0.1$ and $10$ for TeCPSGD. 
\begin{figure}[t]
	\vspace*{-2cm} 
    \begin{minipage}[b]{0.50\linewidth}
    \centering
    \centerline{\includegraphics[width=0.95\linewidth]{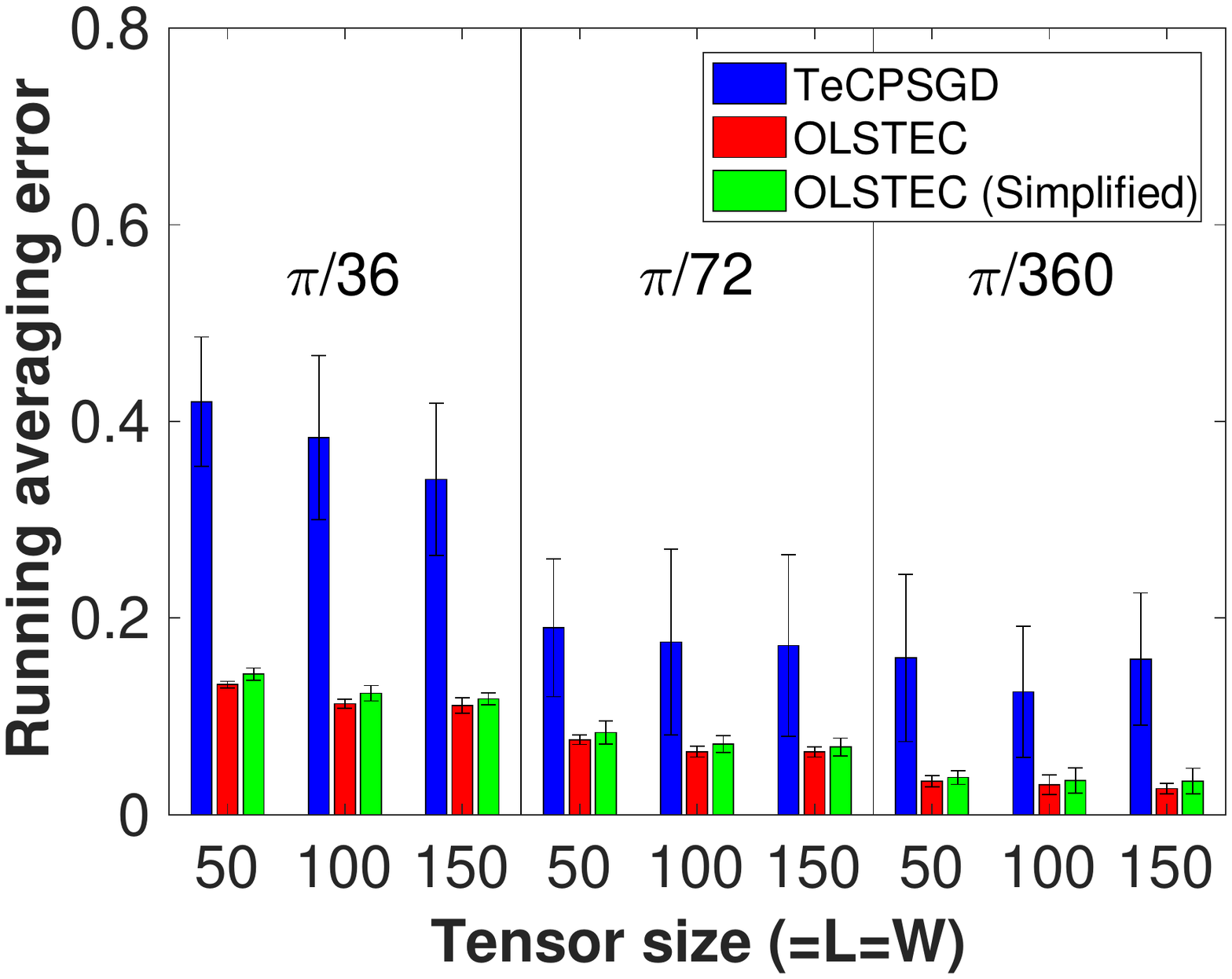}}
    \vspace*{-2cm}
    
    \centerline{\footnotesize (a) $\rho=0.3$}\medskip
    \end{minipage} 
    \begin{minipage}[b]{0.50\linewidth}
    \centering
    \centerline{\includegraphics[width=0.95\linewidth]{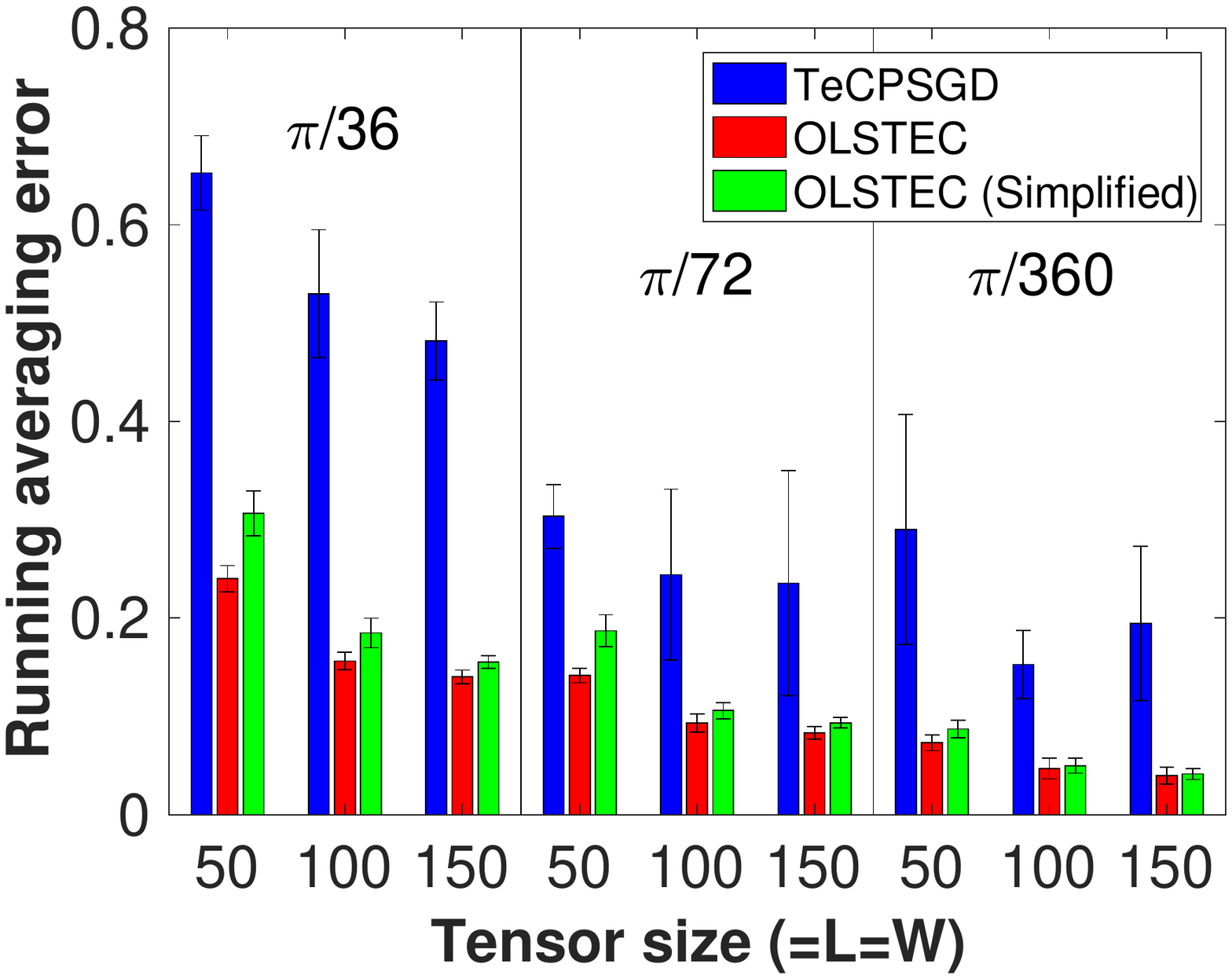}}
    \vspace*{-2cm}
        
    \centerline{\footnotesize (b) $\rho=0.1$}\medskip
    \end{minipage}

    \caption{Running averaging error in synthetic dataset.}
    \label{fig:synthetic_data_ruuningaveraging}
\end{figure}

\begin{table}[t]
\begin{center}
\caption{Processing time [sec] in synthetic dataset.}
\label{tbl:synthetic_data_processing_speed}
\begin{tabular}{c|r|r|r|rr}
\hline
Tensor size& Rank & \hspace*{0.1cm}TeCPSGD\hspace*{0.1cm} & \multicolumn{3}{|c}{OLSTEC}  \\
\cline{4-6}
 $L,W$& $R$\ \ \  &  & \hspace*{0.3cm}Original\hspace*{0.3cm} & Simplified  & (ratio)\ \ \\
\hline\hline
 & 5 & 6.0& 20.6& 17.4& (84.5\%)\\
\cline{2-6}
 & 10 & 6.7& 24.6& 19.7& (80.1\%) \\
\cline{2-6}
 150& 20 & 8.3& 35.8& 24.4& (68.2\%)\\
\cline{2-6}
 & 40 &13.9 & 74.1& 39.5& (53.3\%)\\
\cline{2-6}
 & 60 & 19.2& 131.5& 63.1& (48.0\%)\\
 \hline\hline
& 10 & 64.2& 160.5& 143.4& (89.4\%)\\
\cline{2-6}
 & 20 & 80.9& 209.2& 182.0& (87.0\%) \\
\cline{2-6}
500  & 50 & 157.7& 521.8& 348.6& (66.8\%)\\
\cline{2-6}
 & 100 &492.3 & 1407.5& 865.2& (61.5\%)\\
\cline{2-6}
 & 150 & 776.1& 2817.3& 1676.4& (59.5\%)\\
\hline
\end{tabular}
\end{center}
\end{table}

Figure \ref{fig:synthetic_data_ruuningaveraging} shows the running averaging error for the observation ratios $\rho=\{0.3, 0.1\}$, tensor sizes  $L=W=\{50, 100,150\}$, and angles $\alpha=\{\pi/36, \pi/72, \pi/360\}$, where 10 runs are performed independently. Results show the averages with standard deviations. 
As the figures show, the proposed OLSTEC algorithm shows much lower estimation error, especially when observation ratios are lower. In addition, the standard derivations are also smaller. For that reason, the convergence property of the proposed algorithm is stabler than that of TeCPSGD. Furthermore, comparison of the original OLSTEC algorithm and its simplified version reveals that the simplified OLSTEC algorithm gives similar errors as OLSTEC, especially when the tensor size is large.

Table \ref{tbl:synthetic_data_processing_speed} shows the processing time where the tensor size $L=W$ is $\{150, 500\}$, $T=1000$, and $\rho=0.3$. The rank $R$ is $\{5,10,20,40,60\}$ and $\{10,20,50,100,150\}$ for each tensor size. The time ratio of the simplified OLSTEC algorithm compared with the original algorithm is also shown in brackets. As the figures show, the processing time of TeCPSGD is lower than OLSTEC as expected. The simplified OLSTEC shows much faster than the original OLSTEC, especially when the rank $R$ is larger. Considering the results of the running averaging error together, we conclude that the simplified OLSTEC is preferred to the original OLSTEC for the practical subspace tracking applications. 

%
%
%
\begin{figure}[t]
	\vspace*{-1.5cm} 
    \begin{minipage}[b]{0.333\linewidth}
    \centering
    \centerline{\includegraphics[width=1\linewidth]{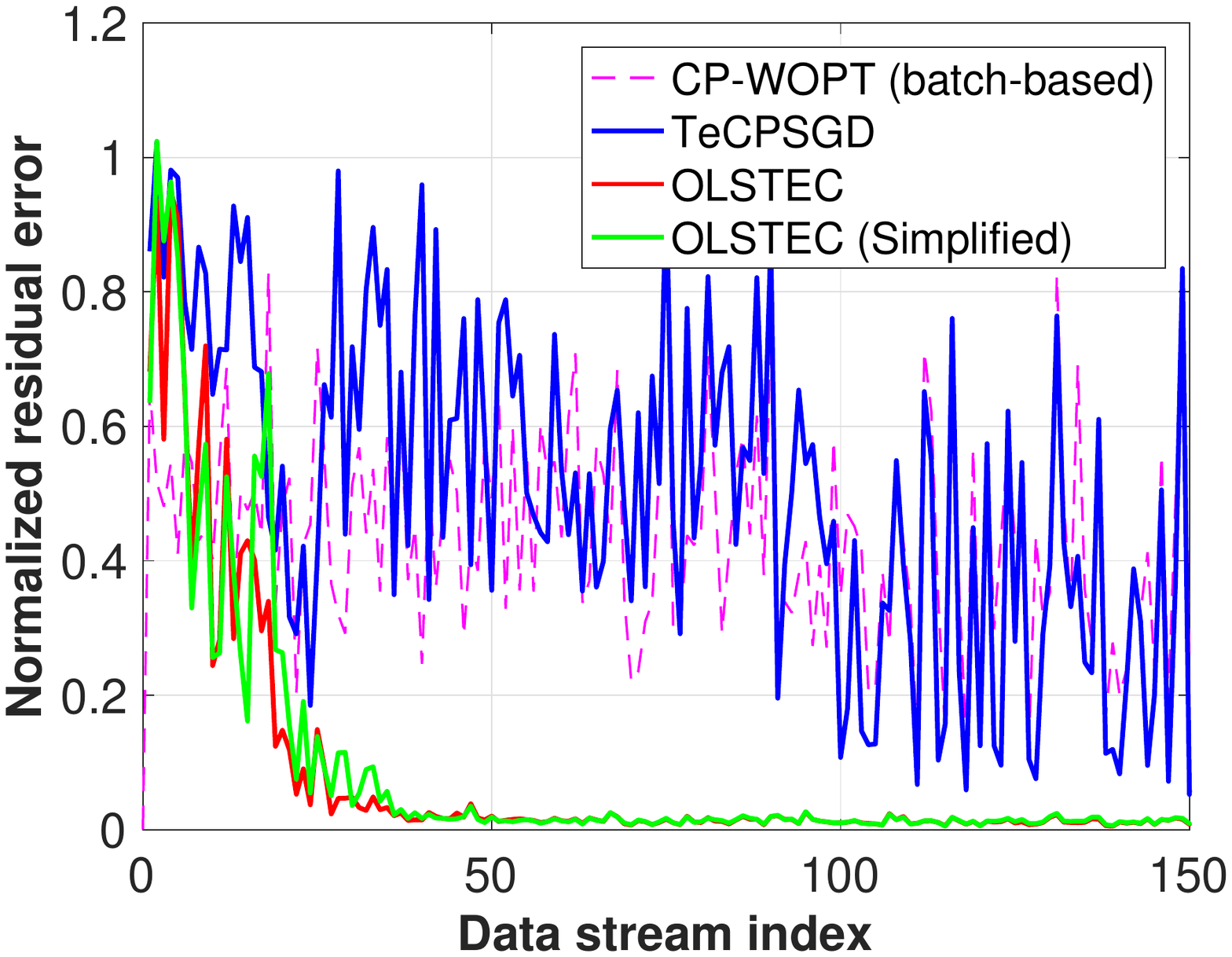}}
    \vspace*{-1.5cm}    
    \centerline{\footnotesize (i) $L=W=50$}\medskip
        
    \end{minipage} 
    \begin{minipage}[b]{0.333\linewidth}
    \centering
    \centerline{\includegraphics[width=1\linewidth]{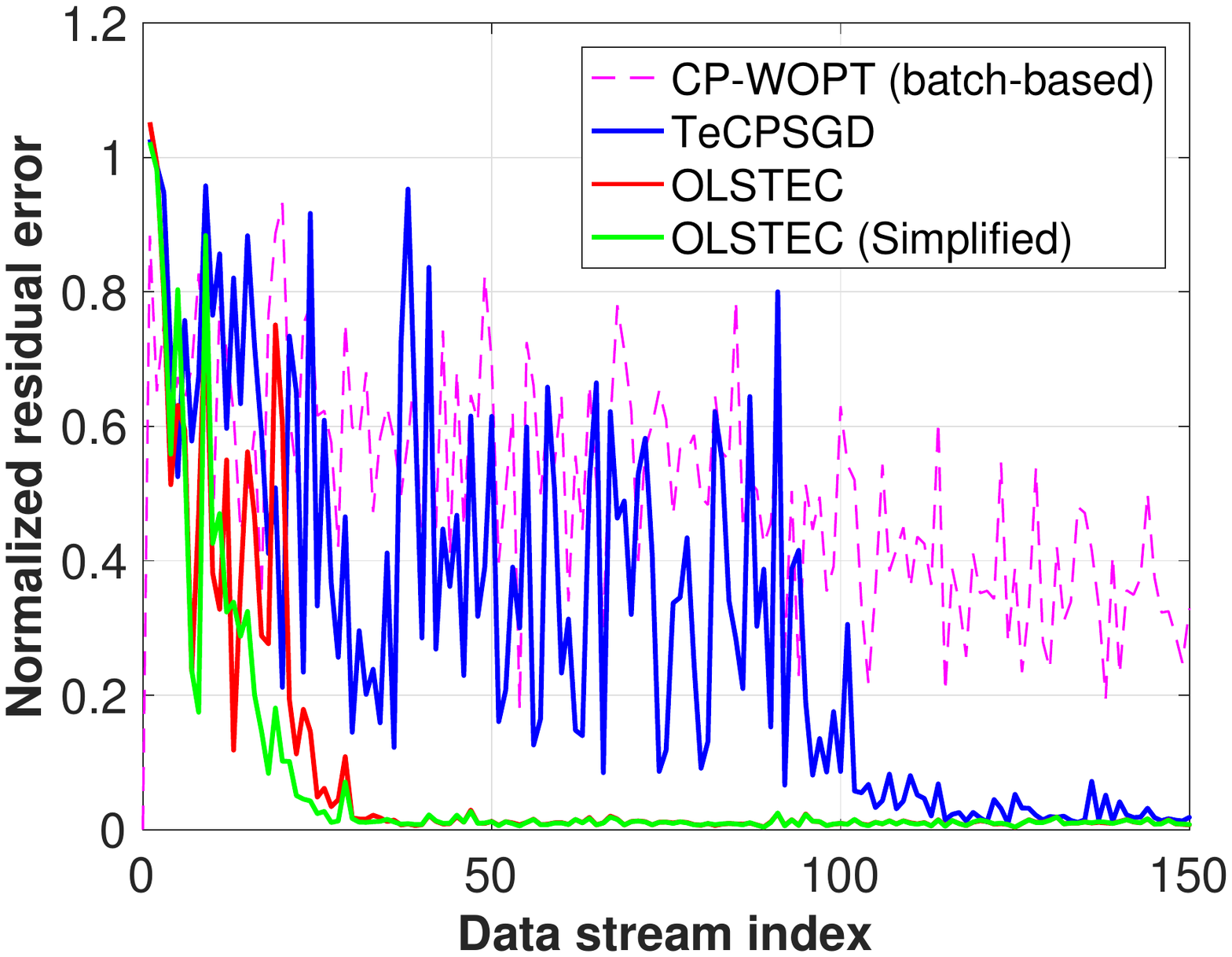}}
    \vspace*{-1.5cm}        

    \centerline{\footnotesize (ii) $L=W=100$}\medskip
    \end{minipage}
    \begin{minipage}[b]{0.333\linewidth}
    \centering
    \centerline{\includegraphics[width=1\linewidth]{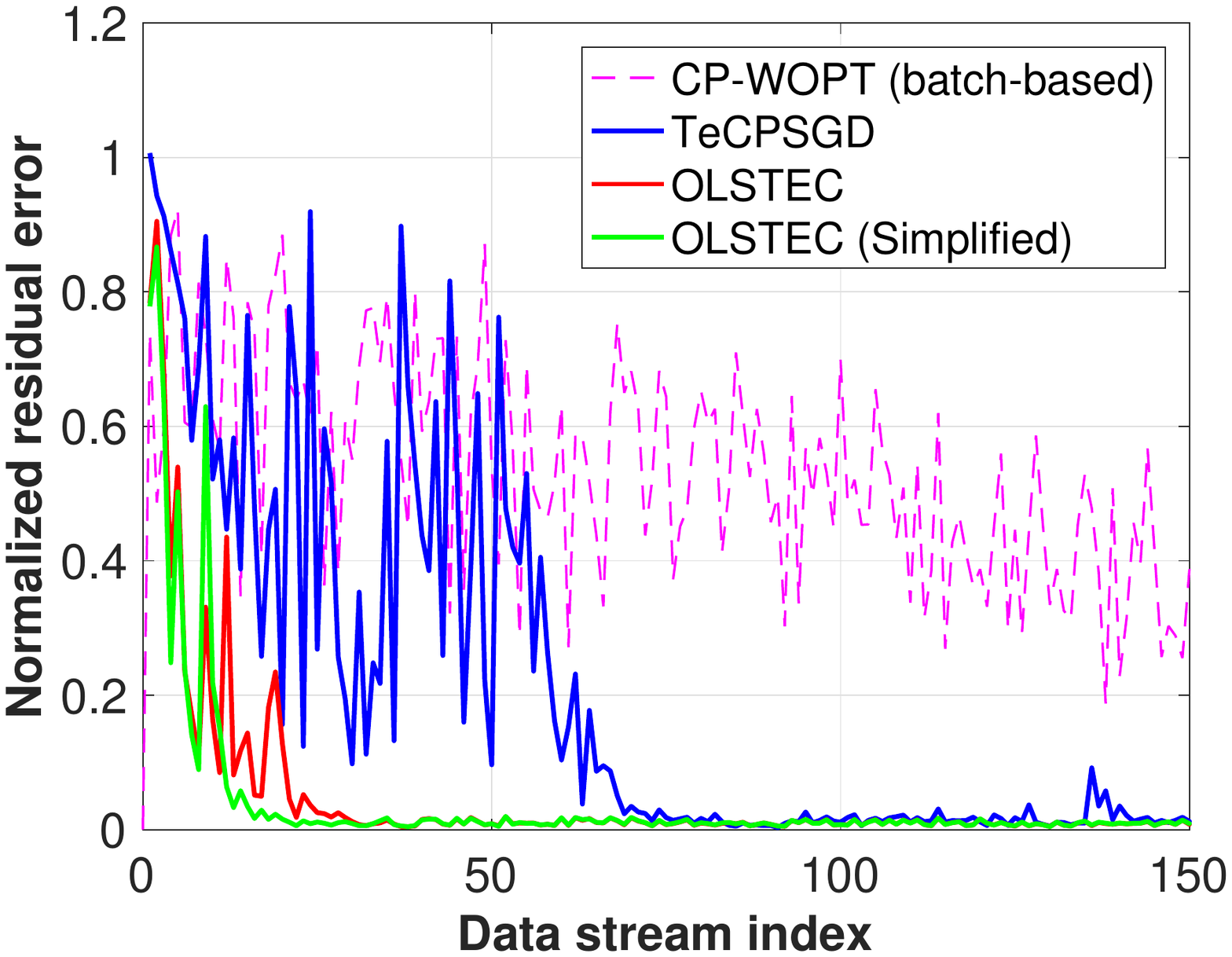}}
    \vspace*{-1.5cm}
        
    \centerline{\footnotesize (iii) $L=W=150$}\medskip
    \end{minipage}

	\centerline{\footnotesize(a) Observation ratio $\rho=0.3$}
	\vspace*{-1.0cm}
	
    \begin{minipage}[b]{0.333\linewidth}
    \centering
    \centerline{\includegraphics[width=1\linewidth]{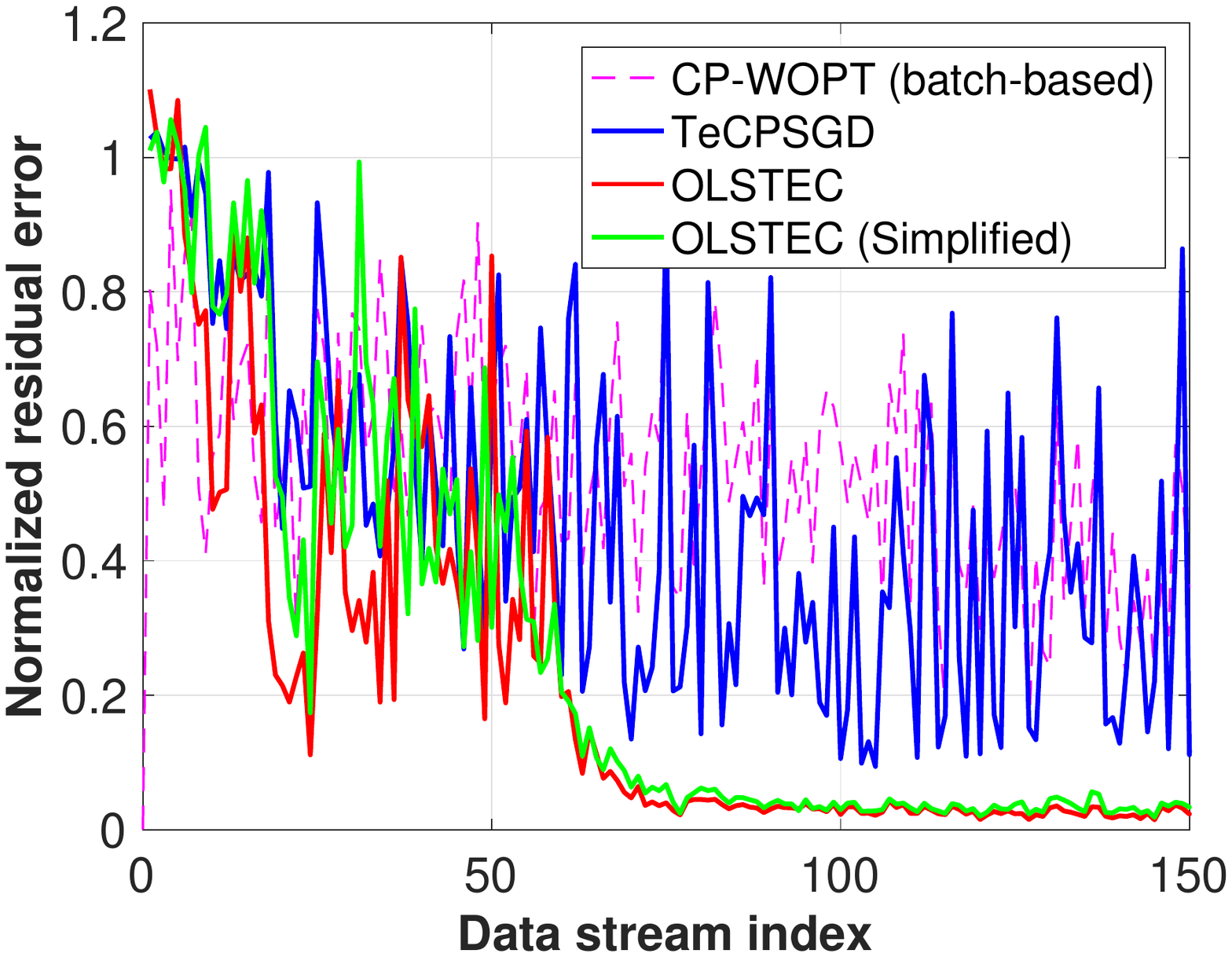}}
    \vspace*{-1.5cm}
        
    \centerline{\footnotesize (i) $L=W=50$}\medskip
    \end{minipage} 
    \begin{minipage}[b]{0.333\linewidth}
    \centering
    \centerline{\includegraphics[width=1\linewidth]{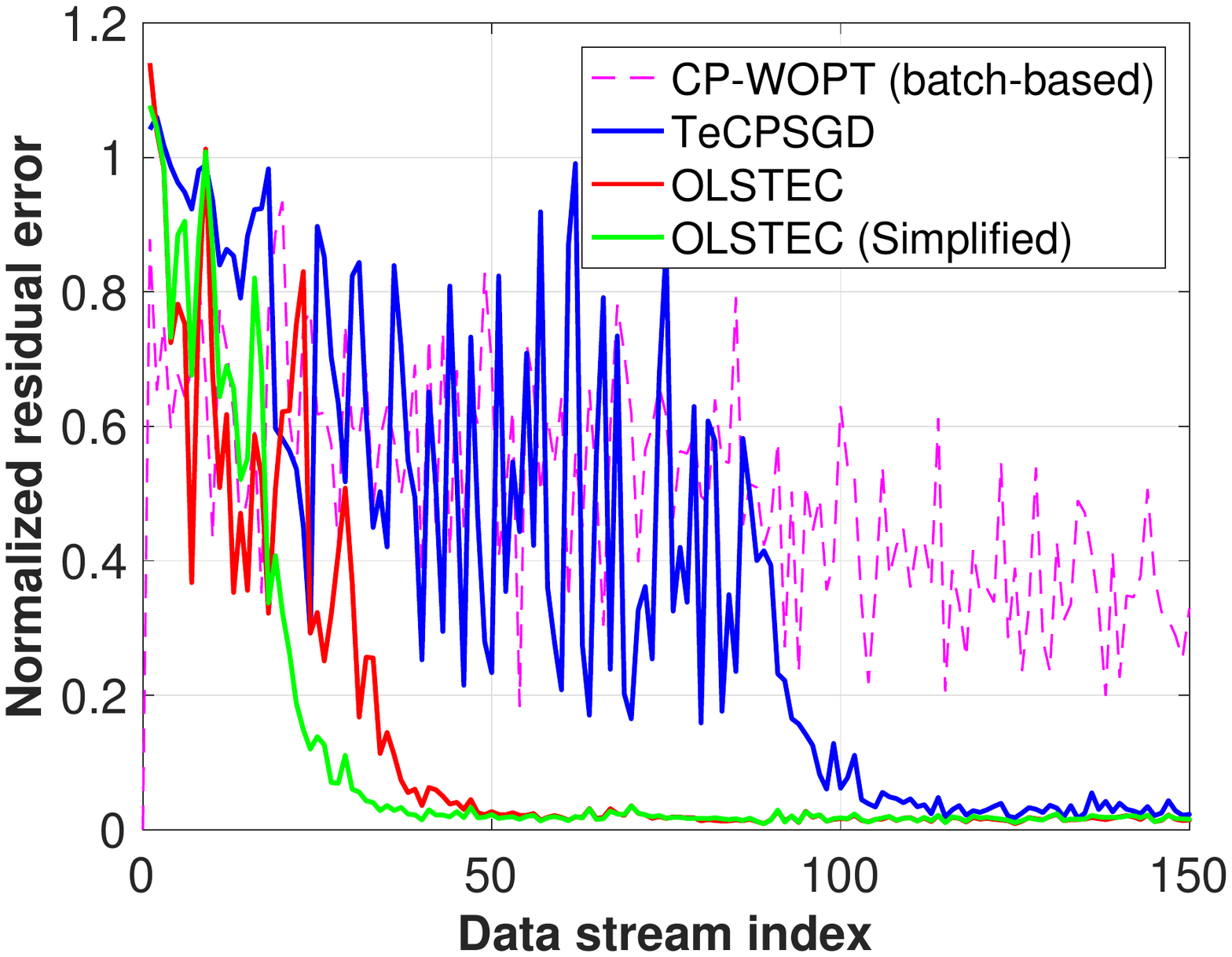}}
    \vspace*{-1.5cm}
        
    \centerline{\footnotesize (ii) $L=W=100$}\medskip
    \end{minipage}
    \begin{minipage}[b]{0.333\linewidth}
    \centering
    \centerline{\includegraphics[width=1\linewidth]{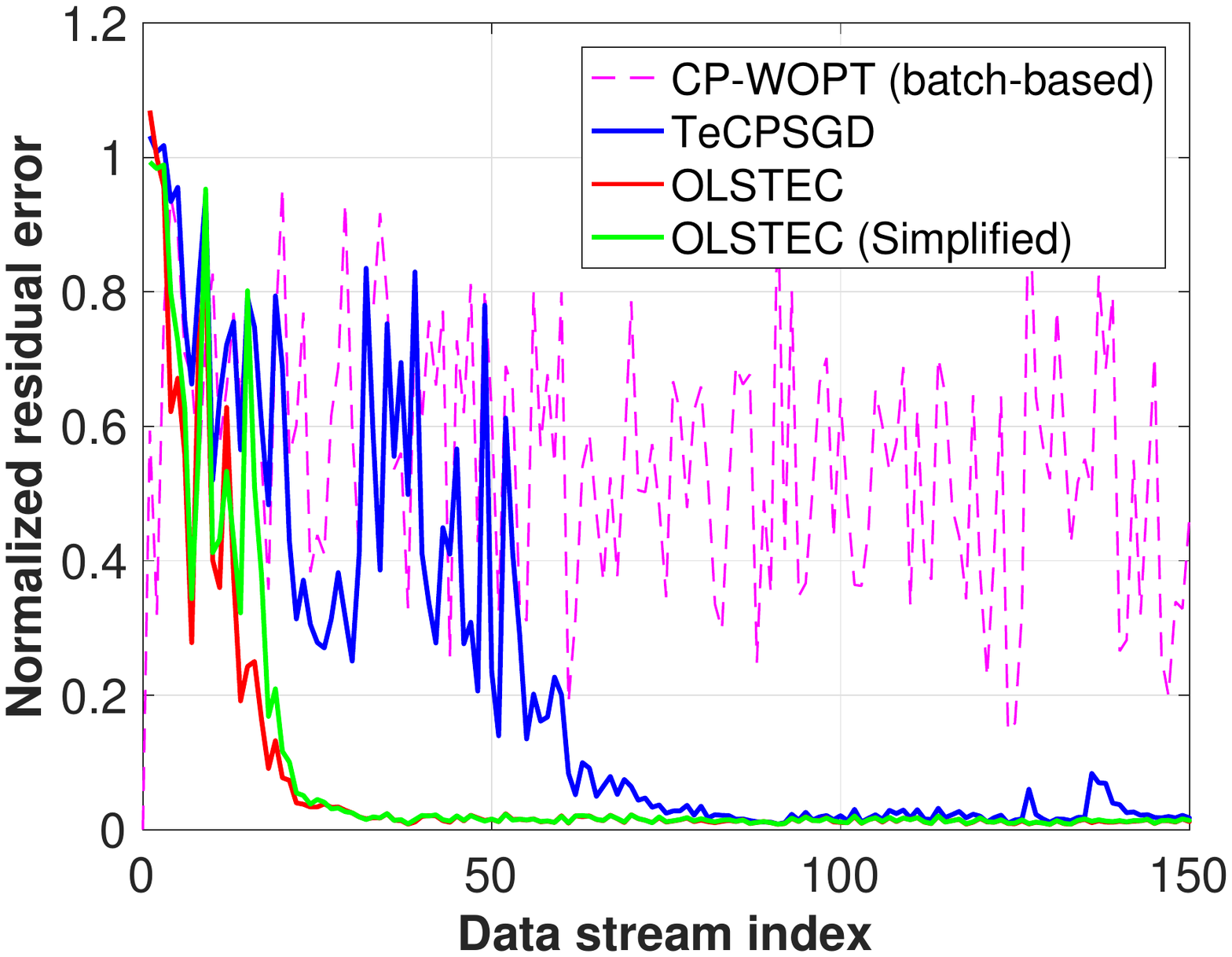}}
    \vspace*{-1.5cm}
        
    \centerline{\footnotesize (iii) $L=W=150$}\medskip
    \end{minipage}
    
	\centerline{\footnotesize(b) Observation ratio $\rho=0.1$}

    \caption{Normalized residual error  in synthetic dataset.}
    \label{fig:synthetic_data_error}
\end{figure}
Figures \ref{fig:synthetic_data_error}(a) and \ref{fig:synthetic_data_error}(b) respectively portray the convergence behaviors of the normalized residual error when the observation ratios are $0.3$ and $0.1$, respectively. 
For this evaluation, we compare OLSTEC with CP-WOPT \cite{Acar_SDM_2010}, the state-of-the-art batch algorithm as reference. The tolerance value of the relative change in function is set to $10^{-9}$. The maximum iterations is $1000$ for CP-WOPT. As expected, CP-WOPT cannot estimate the time-varying subspace because it does not assume that the underlying subspace is time-varying. However, the two proposed algorithms give superior convergence performances than those of CP-WOPT and even TeCPSGD. 
It is also noteworthy that the simplified OLSTEC, surprisingly, outperforms the original OLSTEC in some cases.

\subsection{Network traffic subspace tracking}

We next use real traffic measurements from the GEANT network for validation and evaluation of our proposed algorithm. The GEANT network is a $23$-node network that connects national research and education networks representing $30$ European countries, and does provide internet connectivity to its customers. The dataset consists of four months of traffic matrices that are sampled at the rate of $1/1000$ via {\it Netflow} sampled flow data. The time bin size is $15$ minutes, leading to $672$ sample points for each week's worth of data. We obtain the dataset from the TOTEM project\footnote{\url{https://totem.info.ucl.ac.be/dataset.html}}. In the experiment, we use the tensor data of size $22 \times 22 \times 672$. 
\begin{figure}[t]
\vspace*{-2.5cm}
\begin{center}
\includegraphics[width=8.0cm]{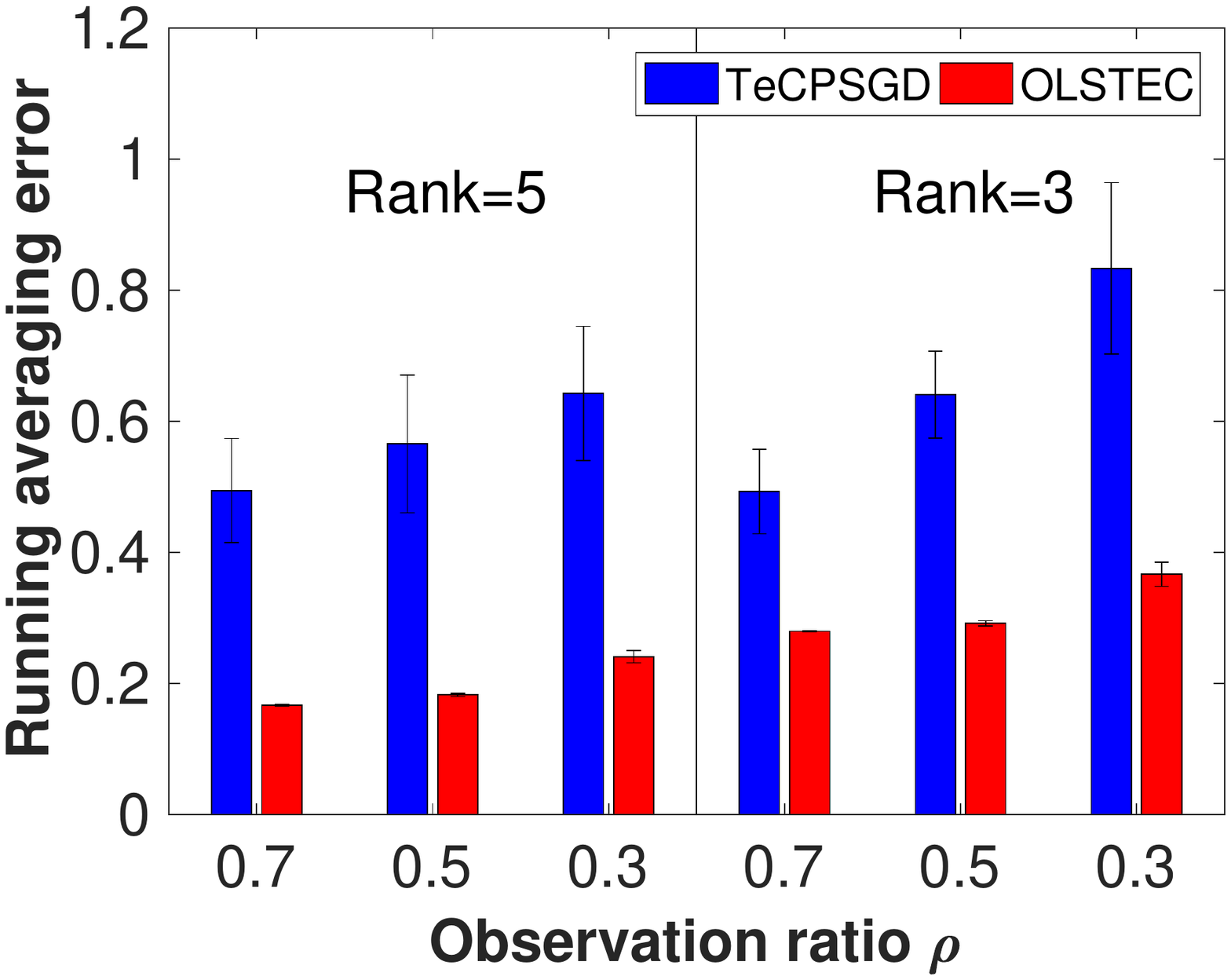}
\vspace*{-2.5cm}

\caption{Running averaging error in the GEANT network traffic.}
\label{fig:ResultGEANT}
\end{center}
\vspace*{-2.5cm}
\begin{minipage}[b]{.5\linewidth}
\centering
\centerline{\includegraphics[width=\linewidth]{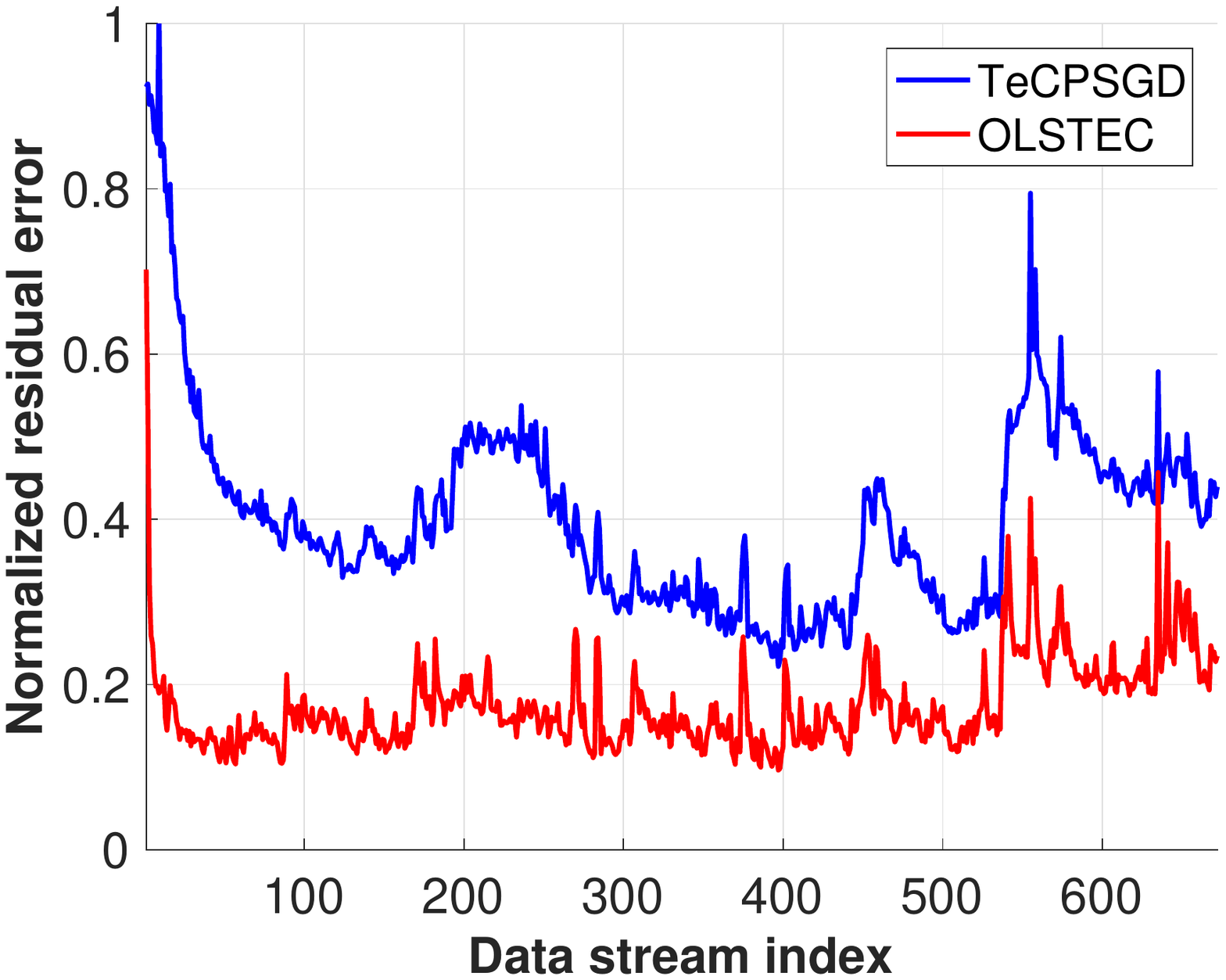}}
\vspace*{-2.5cm}
\centerline{\footnotesize (i) Normalized residual error}
\end{minipage}
\hspace*{0.2cm}
\begin{minipage}[b]{.5\linewidth}
\centering
\centerline{\includegraphics[width=\linewidth]{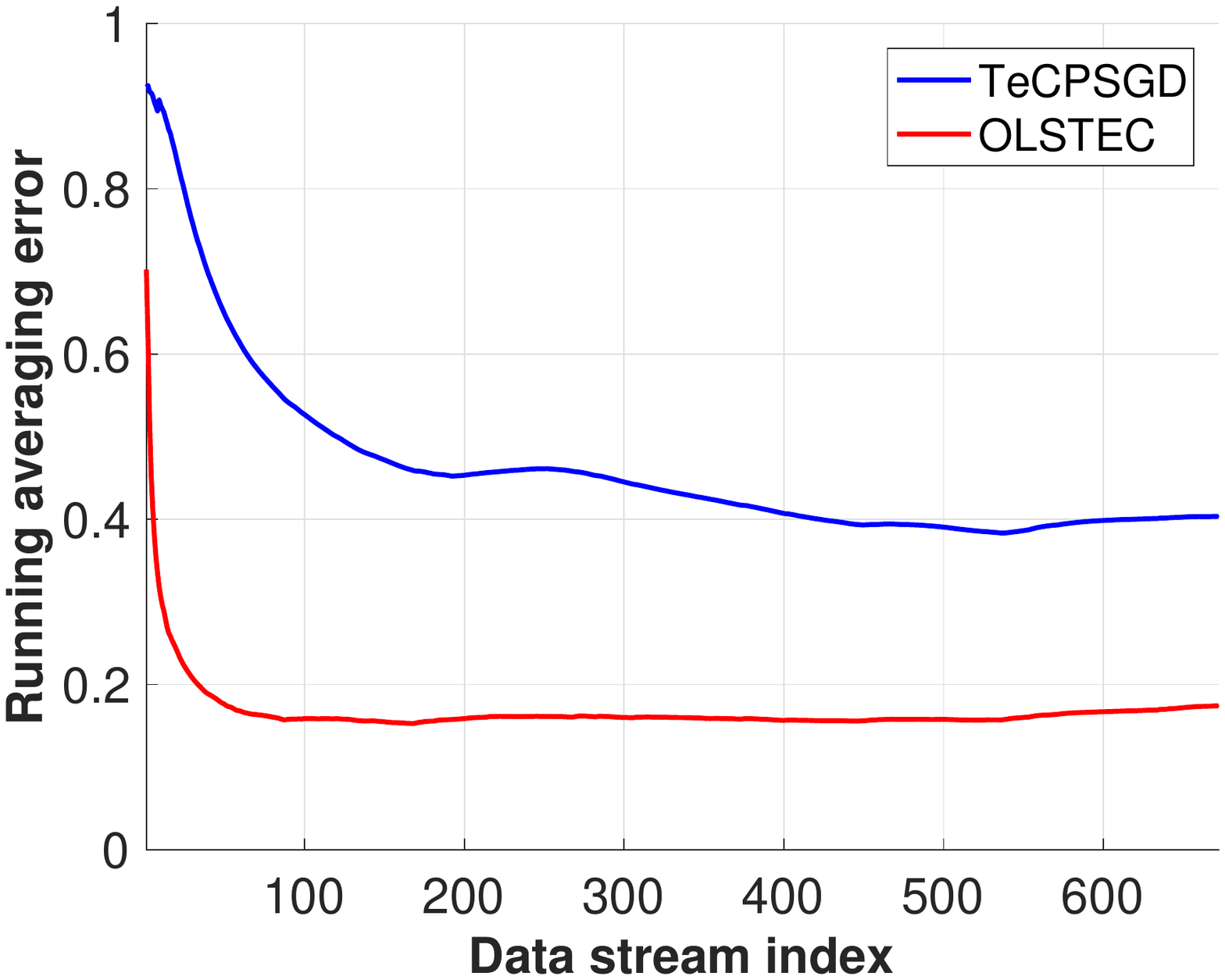}}
\vspace*{-2.5cm}
\centerline{\footnotesize(ii) Running averaging error}
\end{minipage}
\vspace*{0.1cm}

\centerline{\footnotesize(a) Rank $R=5$, Observation ratio $\rho=0.7$}
\label{Fig:ConvErrorGeantNetwork}
\vspace*{-2.cm}

\begin{minipage}[b]{.5\linewidth}
\centering
\centerline{\includegraphics[width=\linewidth]{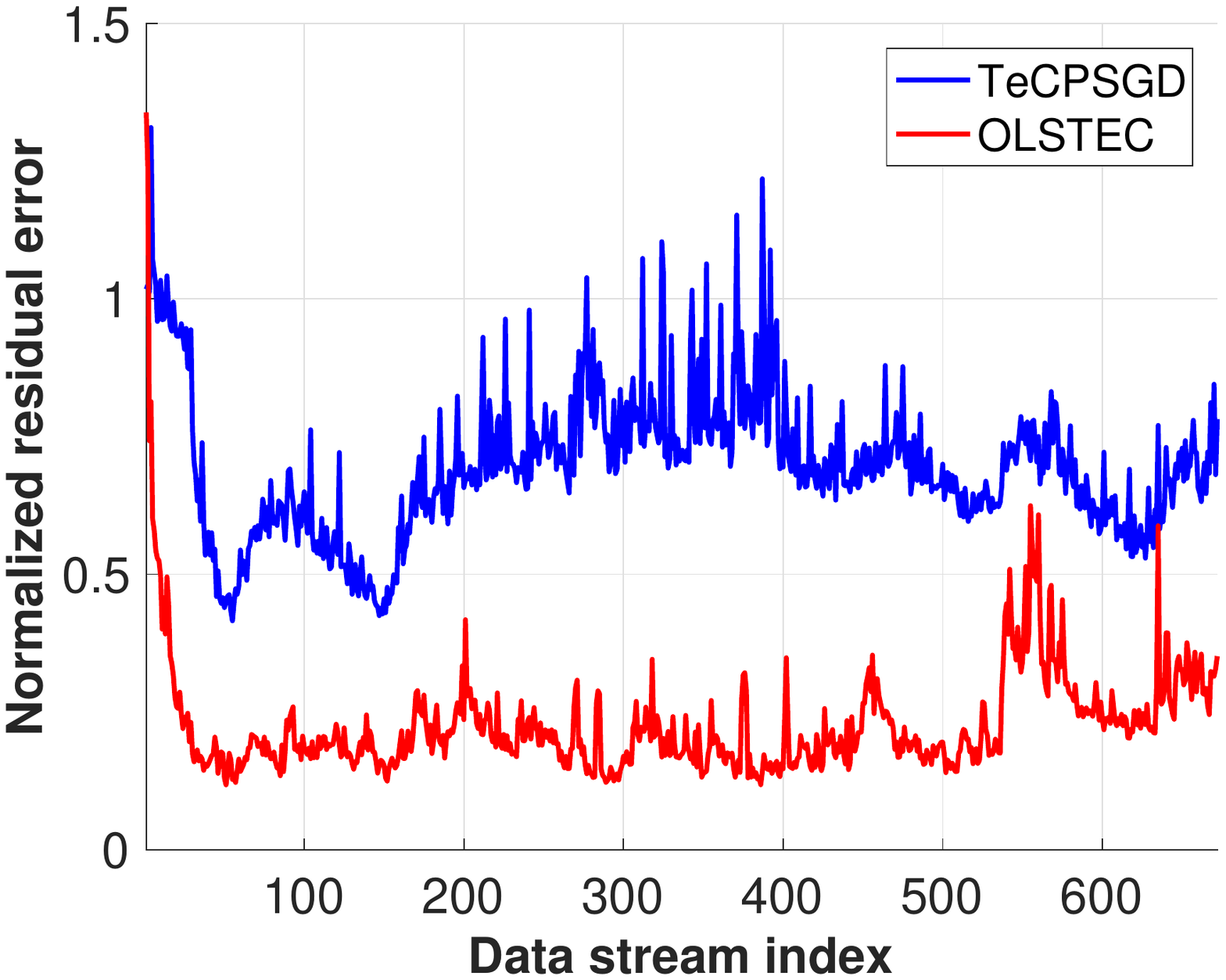}}
\vspace*{-2.5cm}
\centerline{\footnotesize(i) Normalized residual error}
\end{minipage}
\hspace*{0.2cm}
\begin{minipage}[b]{.5\linewidth}
\centering
\centerline{\includegraphics[width=\linewidth]{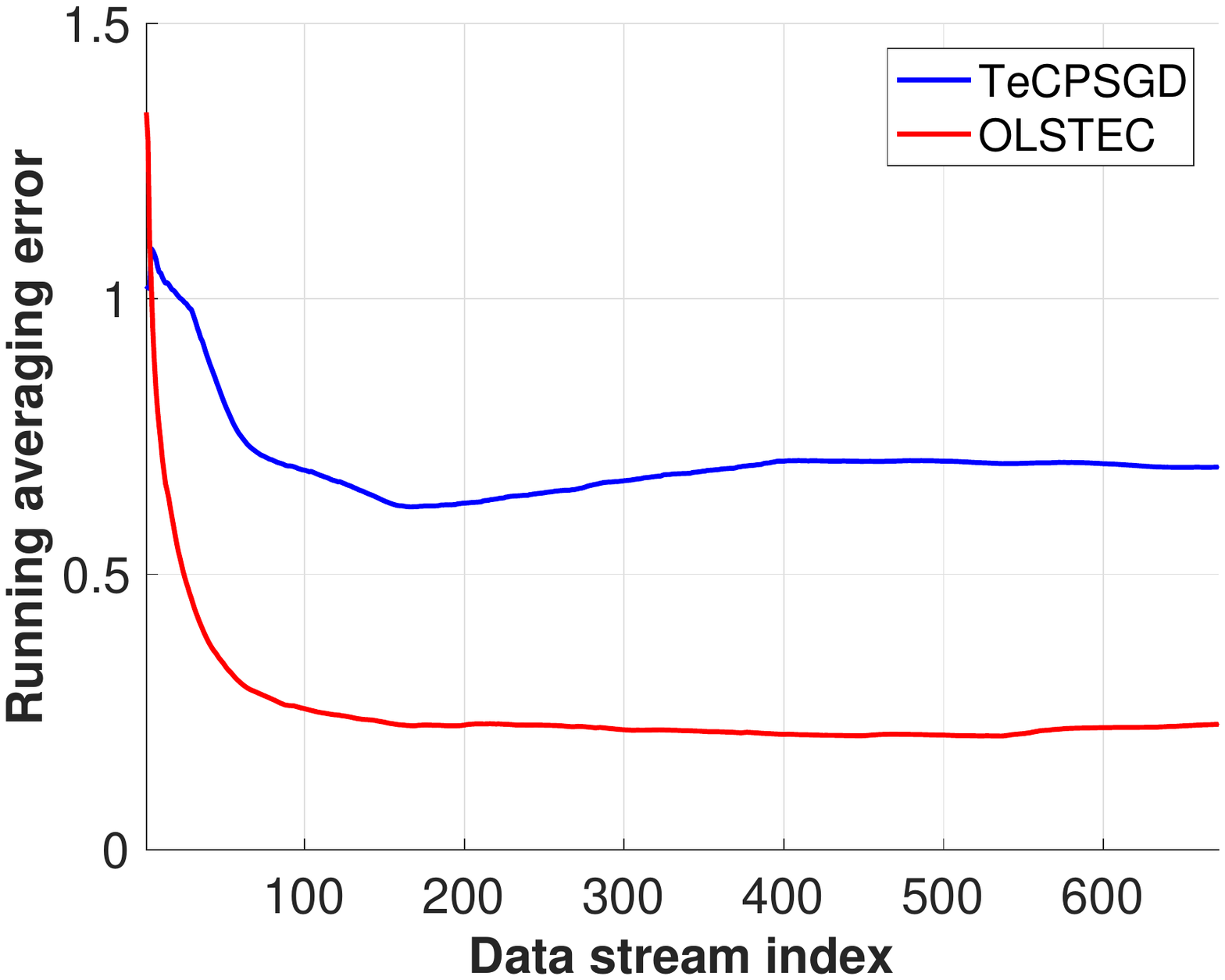}}
\vspace*{-2.5cm}
\centerline{\footnotesize(ii) Running averaging error}
\end{minipage}
\vspace*{0.1cm}

\centerline{\footnotesize(b) Rank $R=5$, Observation ratio $\rho=0.3$}

\caption{Behaviors of errors in the GEANT network traffic.}
\label{Fig:ConvErrorGeantNetwork}
\end{figure}

We evaluate the running averaging error in several settings. The rank $R$ are $\{3,5\}$ and the observation ratios $\rho$ are $\{0.7, 0.5, 0.3\}$. Also, $\lambda$ and $\mu_r$ respectively denoting $0.85$ and $0.1$ are used in the proposed algorithm. $\lambda$, $\mu$ and the stepsize are set respectively to $0.001$, $0.1$ and $0.00001$ for TeCPSGD. These values are obtained from preliminary experiments. The running averaging error of TeCPSGD and our proposed algorithm, OLSTEC, are shown in Figure \ref{fig:ResultGEANT}. These results demonstrate the superior performance of our proposed algorithm. It should be emphasized that the standard deviations of our proposed algorithm are much smaller those of TeCPSGD. These small deviations reflect stabler performance of our proposed algorithm independent of a particular stepsize algorithm. Additionally, the examples of the behavior of the normalized residual error and the running averaging error when $R=5$ and $\rho=\{0.7, 0.3\}$ are presented respectively in Figures \ref{Fig:ConvErrorGeantNetwork}(a) and \ref{Fig:ConvErrorGeantNetwork}(b). As the figures show, the propose algorithm indicates very faster decreases on the errors at the begging of the iteration. Therefore, our proposed algorithm can accommodate the sudden change of the data, thereby providing stabler estimation of the changing subspace than TeCPSGD does.

\clearpage
\subsection{Environmental data subspace tracking}

Next, we evaluate our proposed algorithm for use with environmental datasets using the Metro-UK dataset and the CCDS dataset. 
The Metro-UK dataset is collected from the meteorological office of the UK\footnote{\url{http://www.metoffice.gov.uk/public/weather/climate- historic/}}. It contains monthly measurements of 5 variables in 16 stations across the UK during 1960--2000. Actually, CCDS is the comprehensive climate dataset that is a collection of climate records of North America \cite{Lozano_KDD_2009}. It contains monthly observations of 17 variables such as carbon dioxide and temperature of 125 observation locations during 1990--2001. The sizes of the tensor format in the experiments, $\mathbfcal{Y}$, result respectively in $16 \times 15 \times 492$ and $125 \times 17 \times 156$.  

\begin{figure}[htbp]
\vspace*{-1.5cm}
\begin{minipage}[b]{.32\linewidth}
\centering
\centerline{\includegraphics[width=\linewidth]{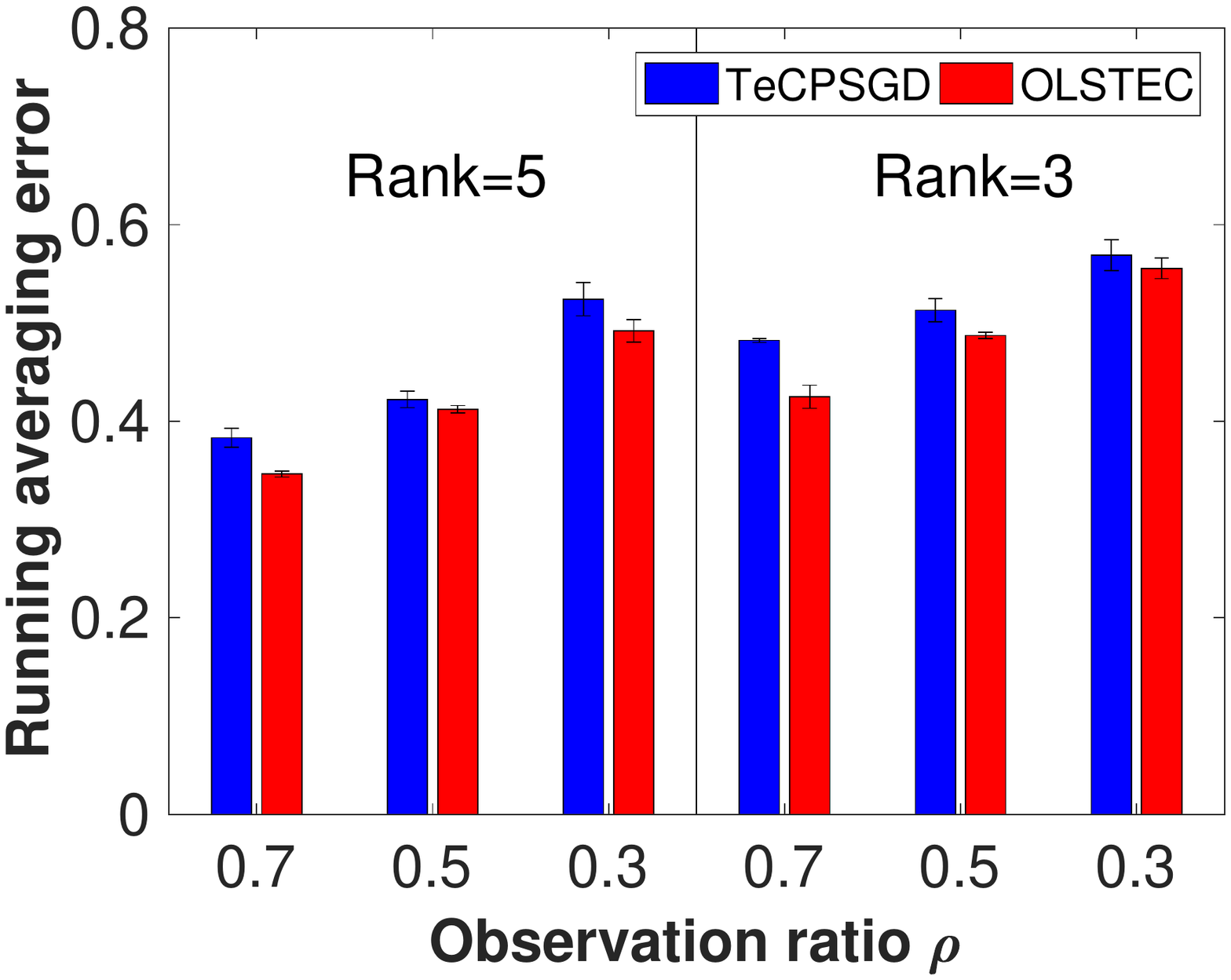}}
\vspace*{-1.5cm}
\centerline{\footnotesize(a)Running averaging error}
\end{minipage}
\hspace*{0.05cm}
\begin{minipage}[b]{.32\linewidth}
\centering
\centerline{\includegraphics[width=\linewidth]{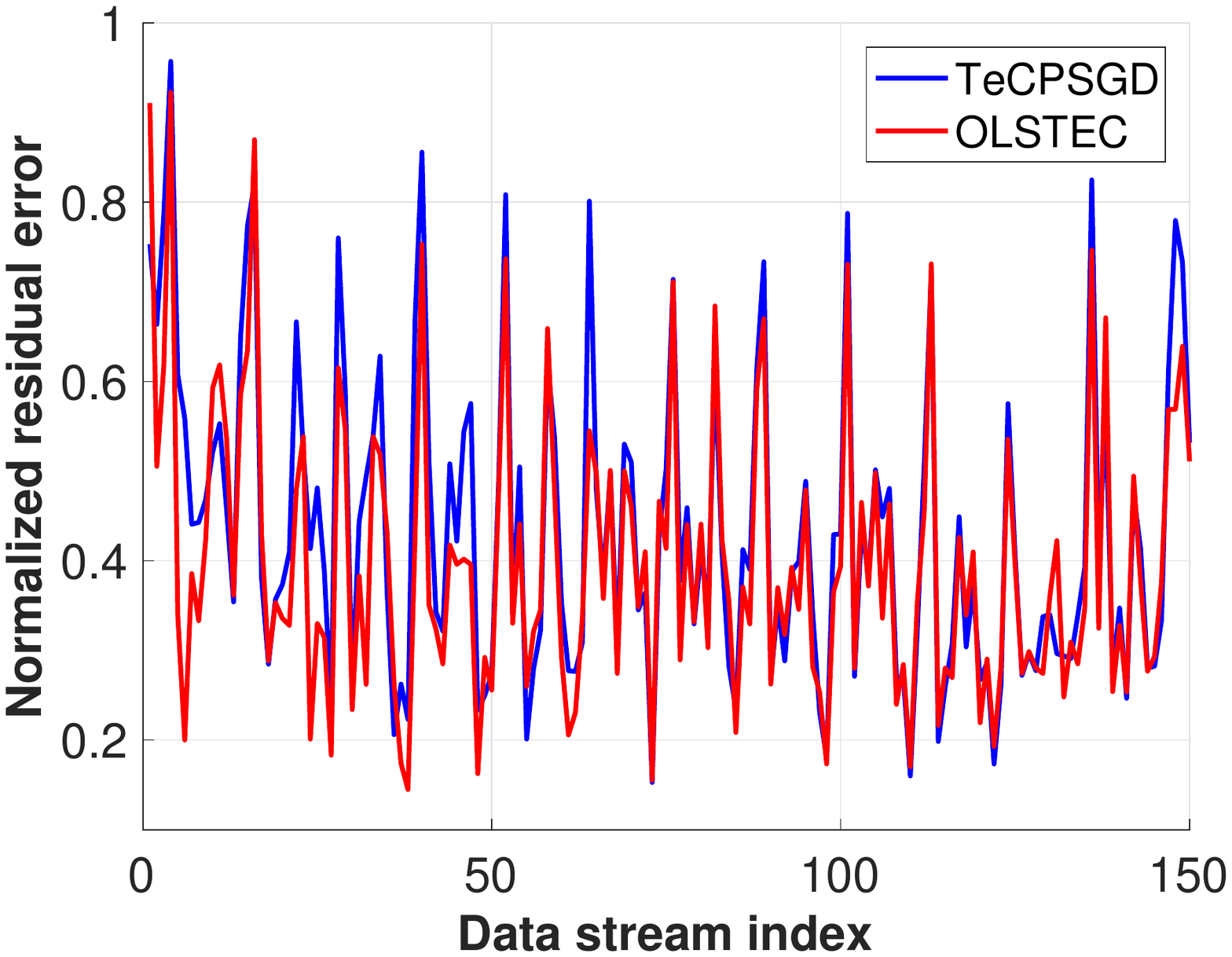}}
\vspace*{-1.5cm}
\centerline{\footnotesize(b)Normalized residual error}
\end{minipage}
\hspace*{0.05cm}
\begin{minipage}[b]{.32\linewidth}
\centering
\centerline{\includegraphics[width=\linewidth]{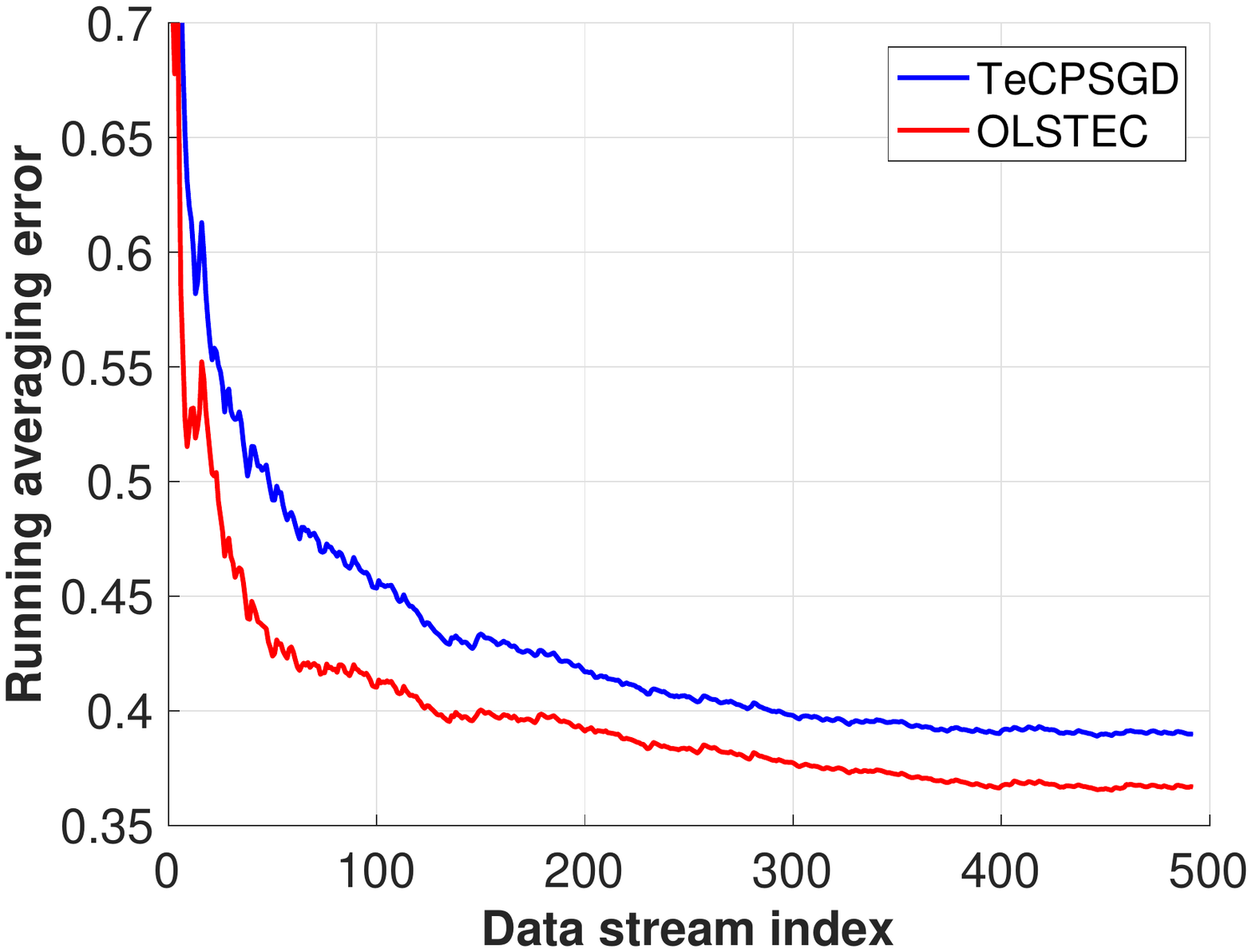}}
\vspace*{-1.5cm}
\centerline{\footnotesize(c)Running averaging error}
\end{minipage}
\caption{Behavior of normalized residual error and running averaging error in UK-metro.}
\label{Fig:ConvErrorUKMetro}
\vspace*{-1.5cm}
%
%
%
\vspace*{0.5cm}
\begin{minipage}[b]{.32\linewidth}
\centering
\centerline{\includegraphics[width=\linewidth]{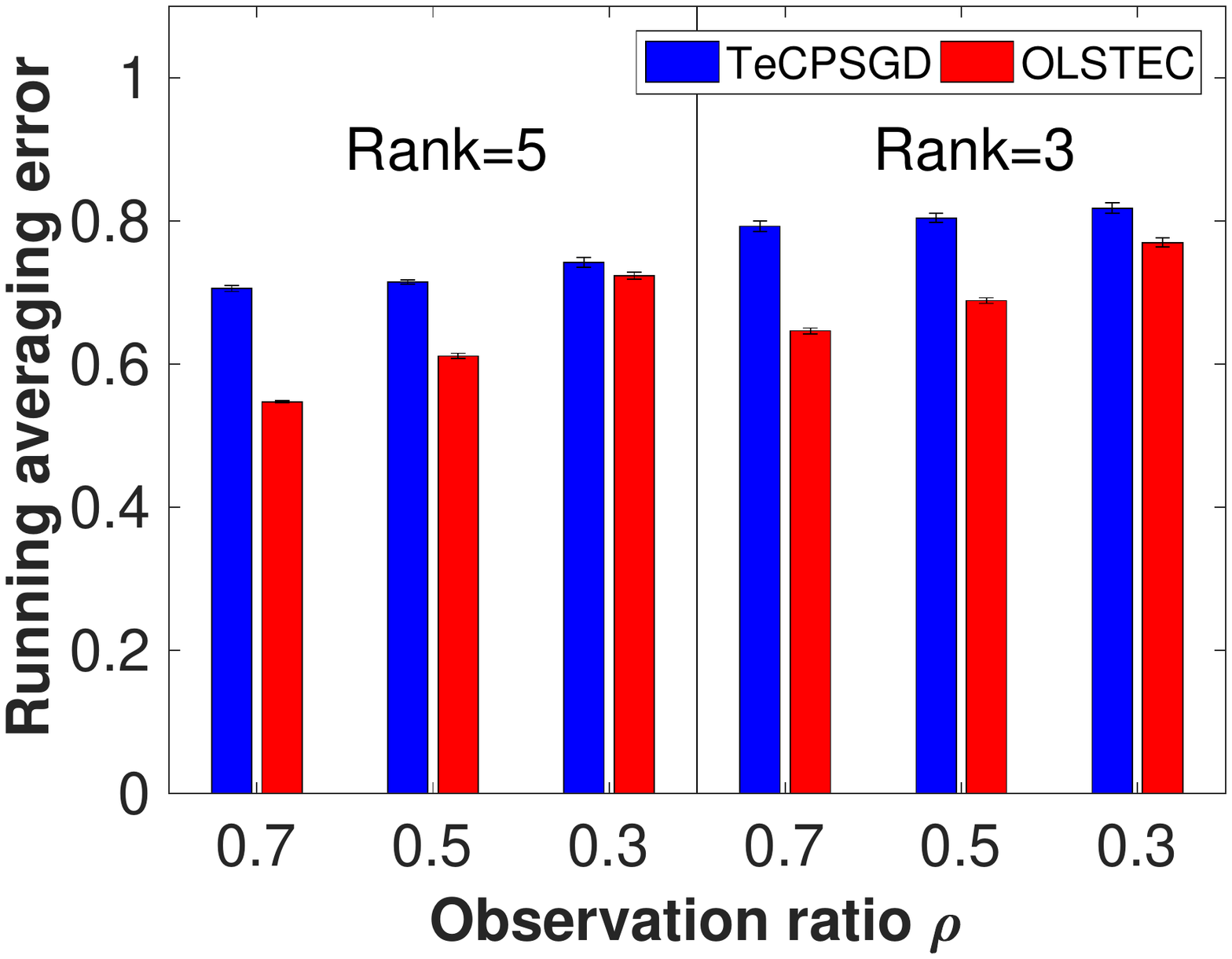}}
\vspace*{-1.5cm}
\centerline{\footnotesize(a)Running averaging error}\medskip
\end{minipage}
\hspace*{0.05cm}
\begin{minipage}[b]{.32\linewidth}
\centering
\centerline{\includegraphics[width=\linewidth]{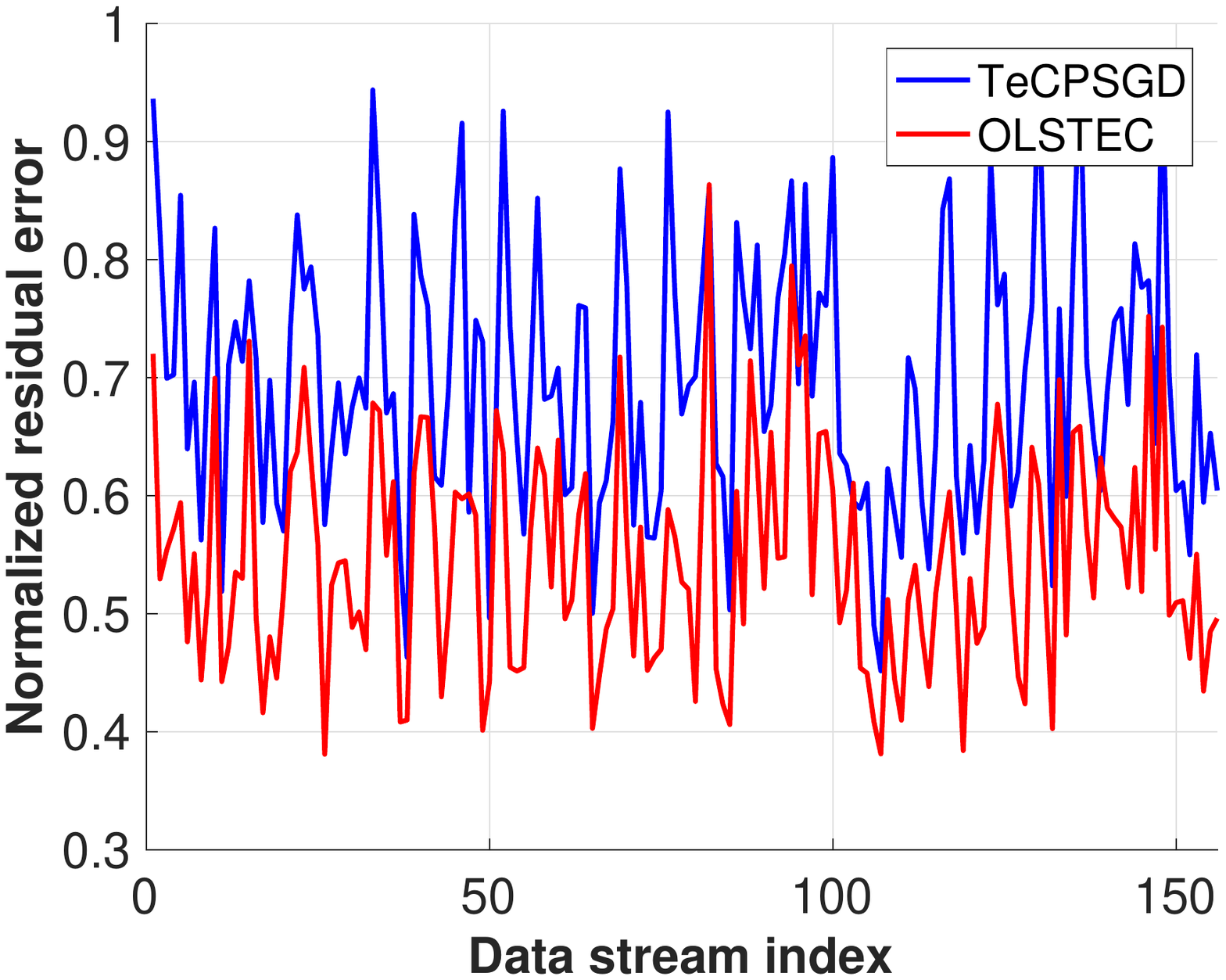}}
\vspace*{-1.5cm}
\centerline{\footnotesize(b)Normalized residual error}\medskip
\end{minipage}
\hspace*{0.05cm}
\begin{minipage}[b]{.32\linewidth}
\centering
\centerline{\includegraphics[width=\linewidth]{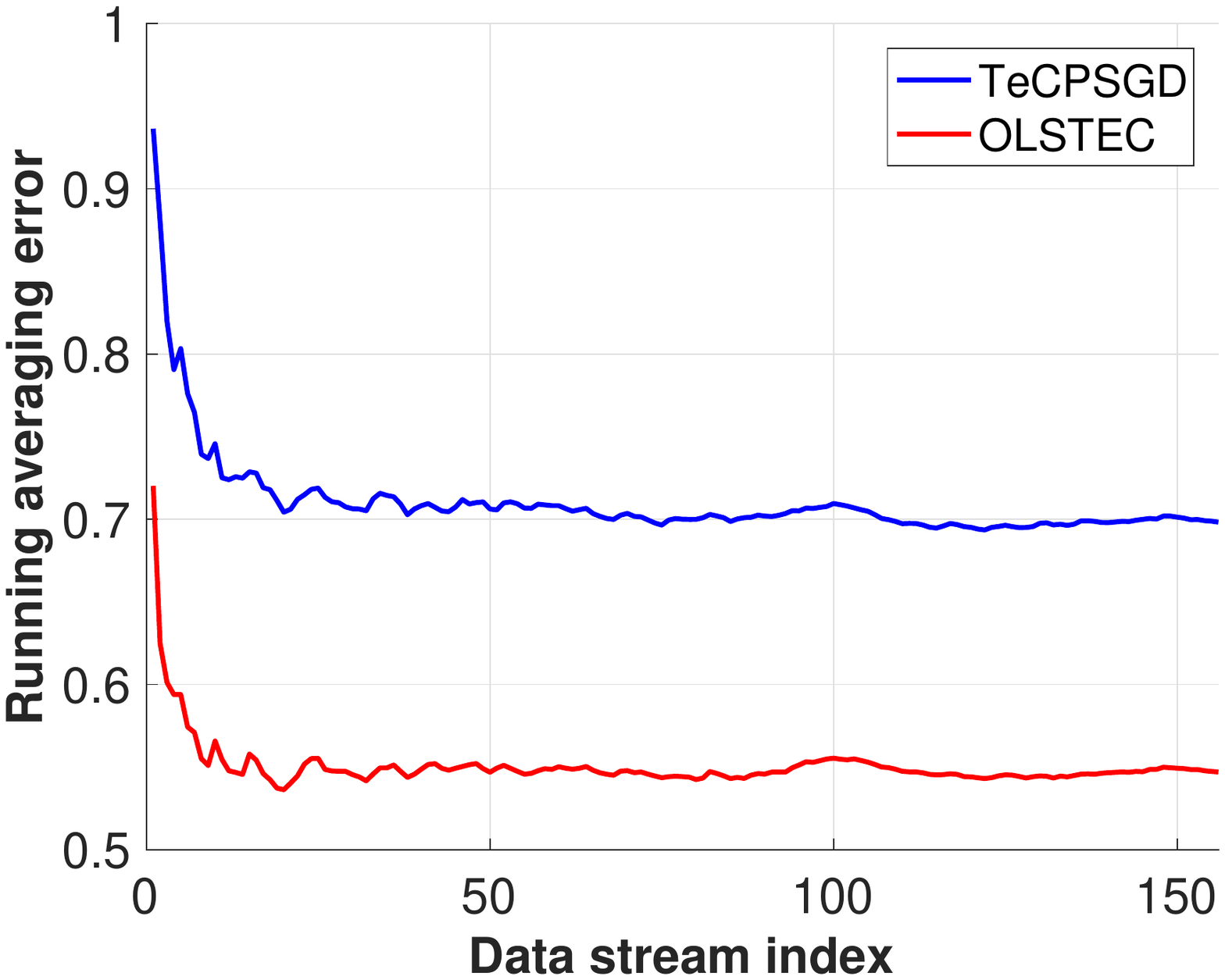}}
\vspace*{-1.5cm}
\centerline{\footnotesize(c)Running averaging error}\medskip
\end{minipage}
\caption{Behavior of normalized residual error and running averaging error in CCDS. }
\label{Fig:ConvErrorCCDS}
\end{figure}

Figures \ref{Fig:ConvErrorUKMetro} and \ref{Fig:ConvErrorCCDS} respectively present for the Metro-UK dataset and the CCDS dataset. 
$\lambda=0.3$ and $\mu_r=0.1$ are used, respectively, in the proposed algorithm. $\lambda$, $\mu$ and the stepsize are set, respectively, to $0.001$, $0.1$ and $10$ for TeCPSGD. These values are also obtained from preliminary experiments. Panels (a) of both figures show the summary of the running averaging error when $\rho=\{0.7, 0.5, 0.3\}$ and the rank $R=5$. These results give the superior performance of our proposed algorithm. In addition, panels (b) and (c) of the respective figures show the change of the normalized residual error every iteration and the running averaging error, respectively. Panel (b) shows only the partial results up to $t=150$ to distinguish the lines. From these results, it is apparent that our proposed algorithm gives better performance than TeCPSGD does. 

\subsection{Video background subtraction tracking}

\begin{figure}[htbp]
\vspace*{-2.5cm}
\begin{minipage}[b]{.5\linewidth}
\centering
\includegraphics[width=1.1\linewidth]{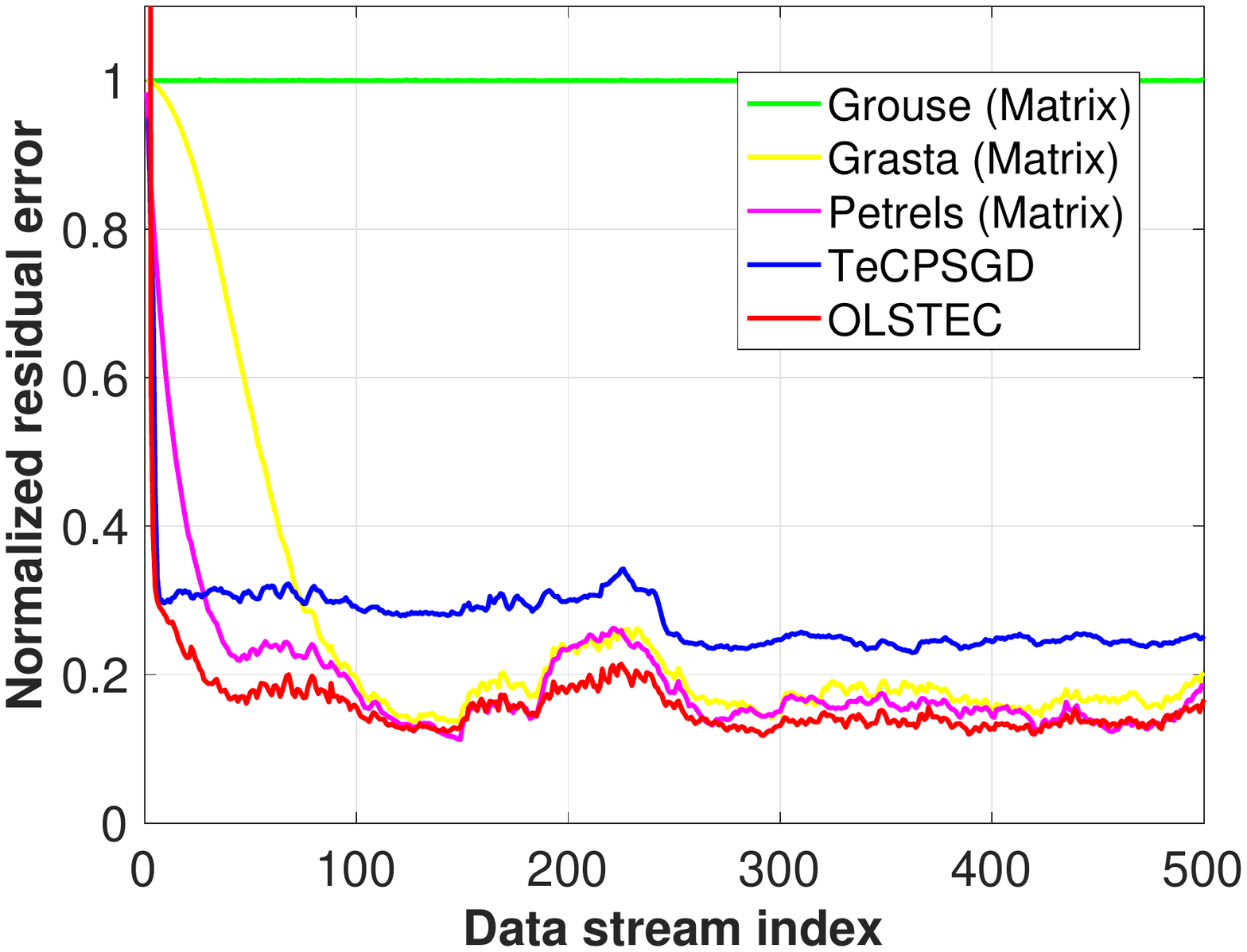}
\vspace*{-3.2cm}

\centerline{\footnotesize (a) Stationary background}
\end{minipage}
\hspace*{0.2cm}
\begin{minipage}[b]{.5\linewidth}
\centering
\includegraphics[width=0.94\linewidth]{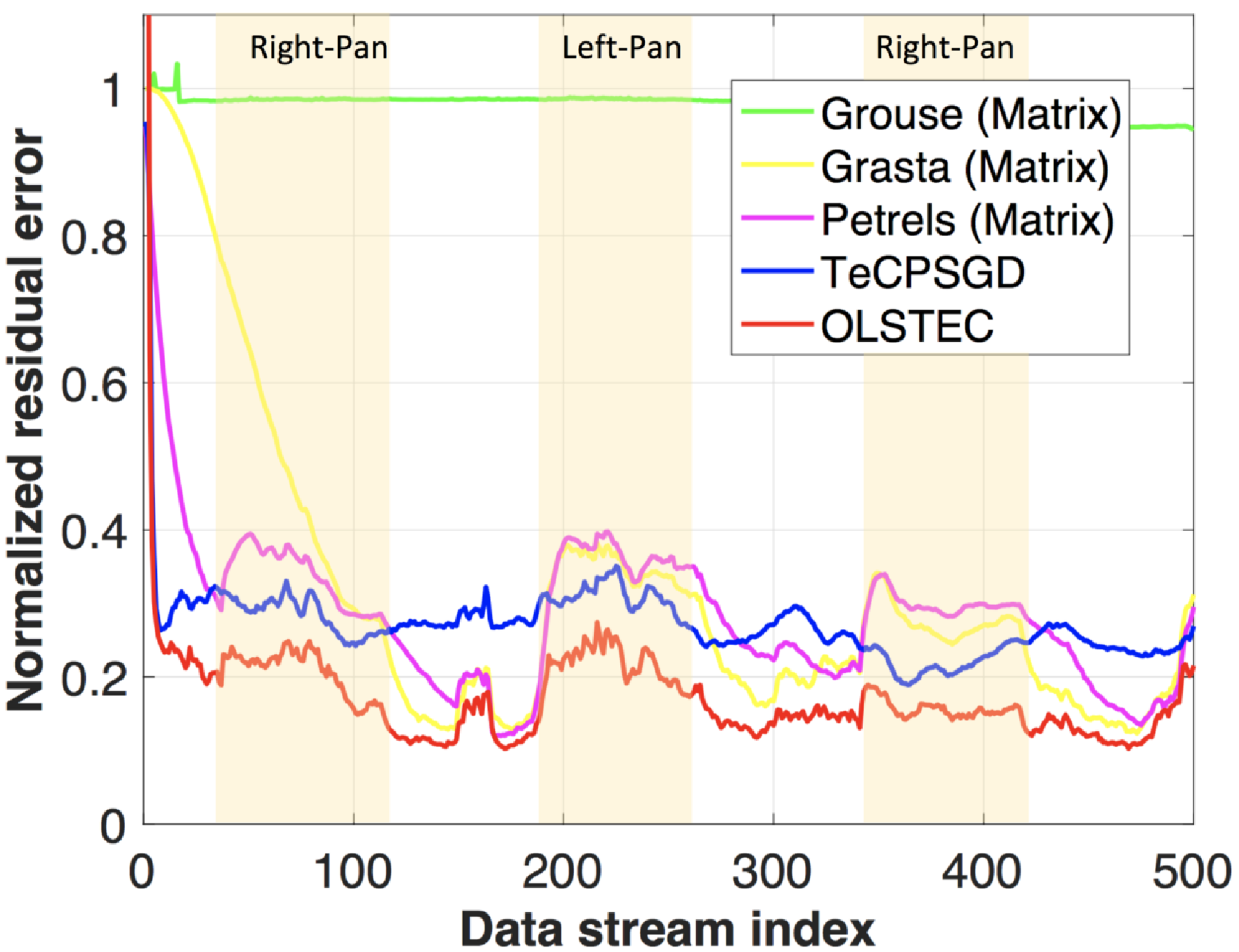}
\centerline{\footnotesize (b) Dynamic background}
\end{minipage}

\caption{The normalized estimation error in surveillance video Airport Hall.}
\label{fig:NRS_Video}
\end{figure}

\begin{figure}[t]
\vspace*{-2cm}
\begin{tabular}{ccc}
\begin{minipage}[b]{.32\hsize}
\centering
\centerline{\includegraphics[width=1\hsize]{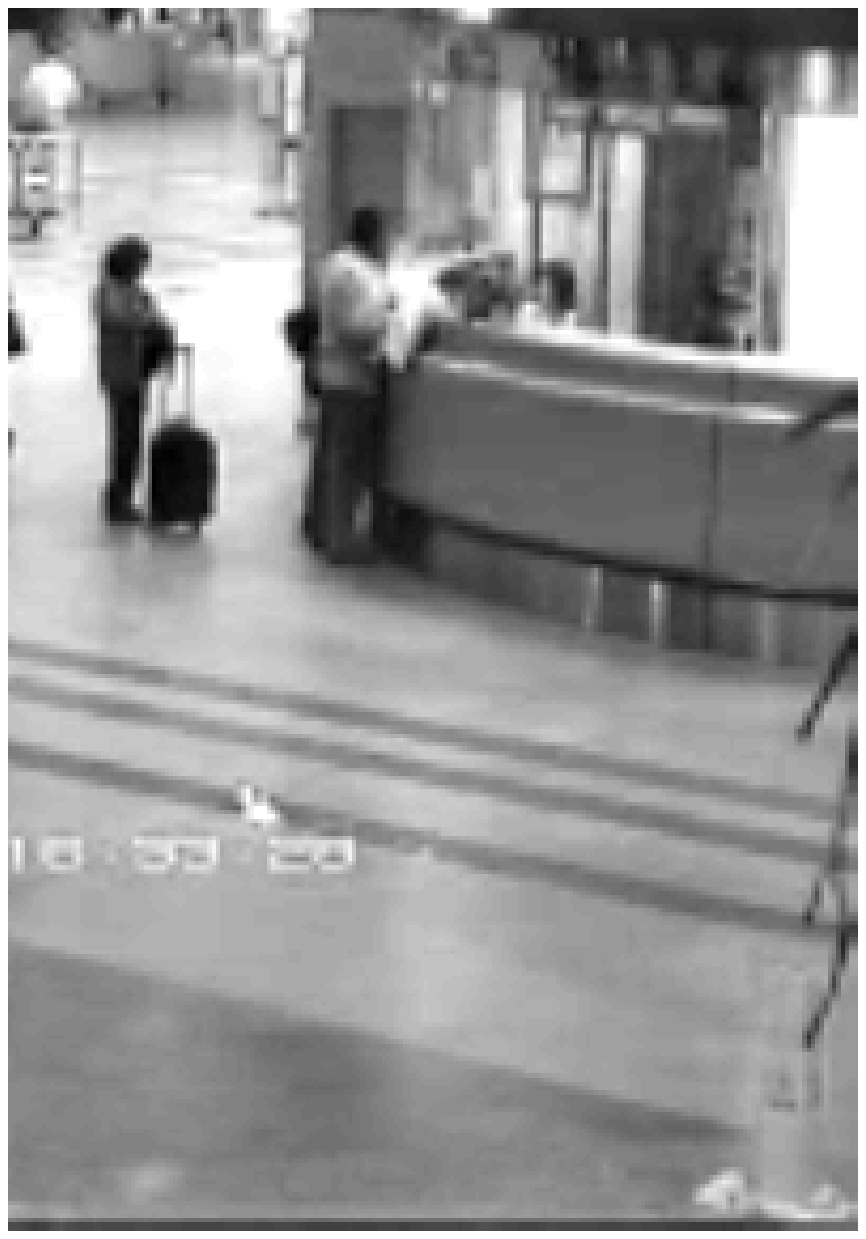}}
\vspace*{-2cm}
\centerline{\footnotesize (i) Original}\medskip
\end{minipage}
&%
\begin{minipage}[b]{.32\hsize}
\centerline{\includegraphics[width=1\hsize]{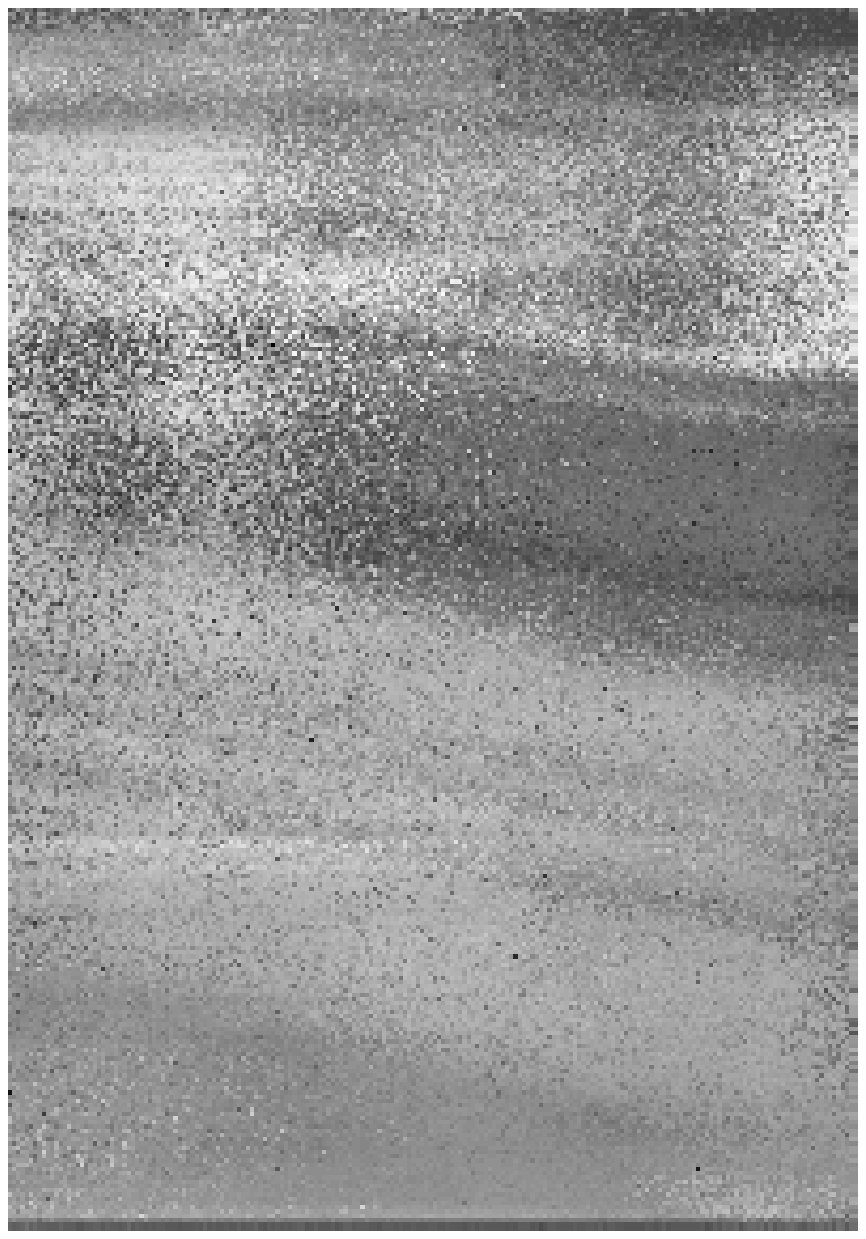}}
\vspace*{-2cm}
\centerline{\footnotesize (ii) GROUSE}\medskip
\end{minipage}
&%
\begin{minipage}[b]{0.32\hsize}
\centerline{\includegraphics[width=1\hsize]{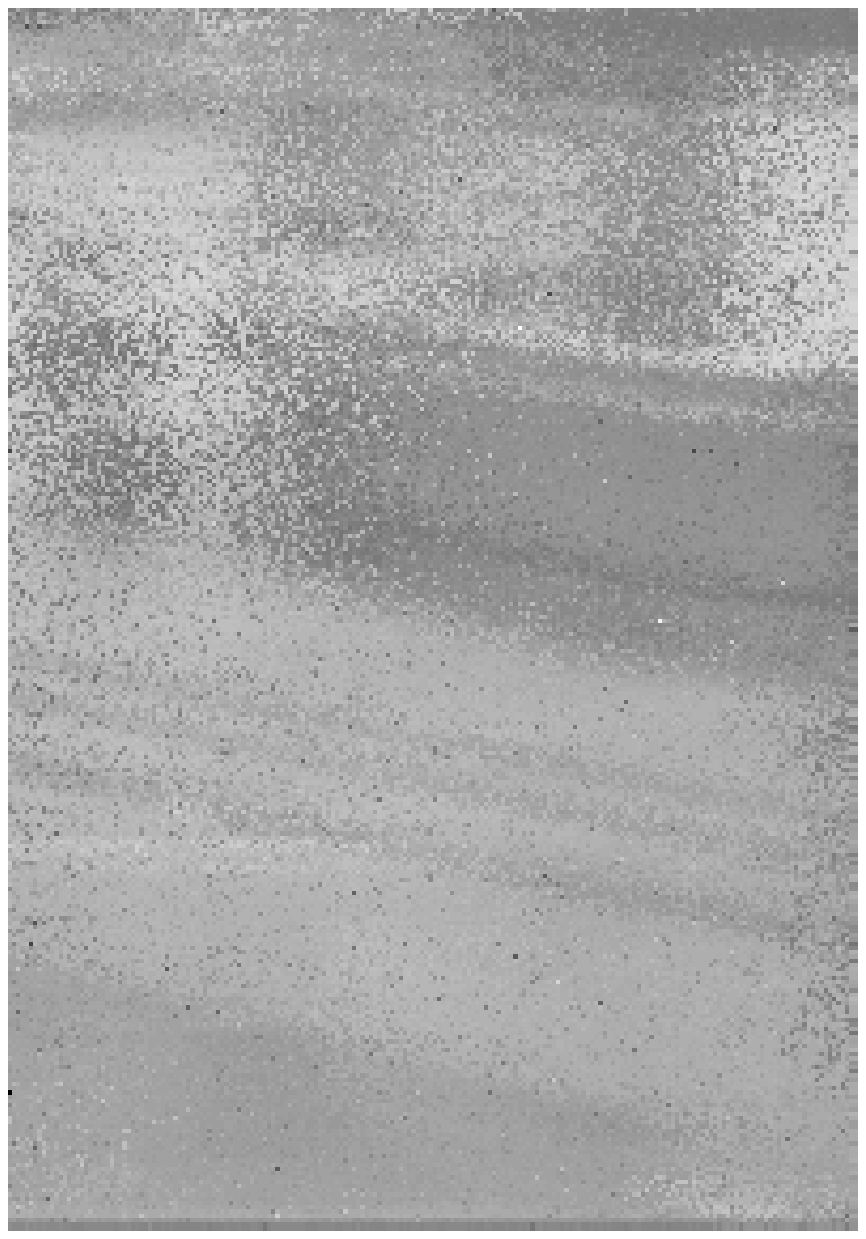}}
\vspace*{-2cm}
\centerline{\footnotesize (iii) GRASTA}\medskip
\end{minipage} \\
\begin{minipage}[b]{.32\hsize}
\centerline{\includegraphics[width=1\hsize]{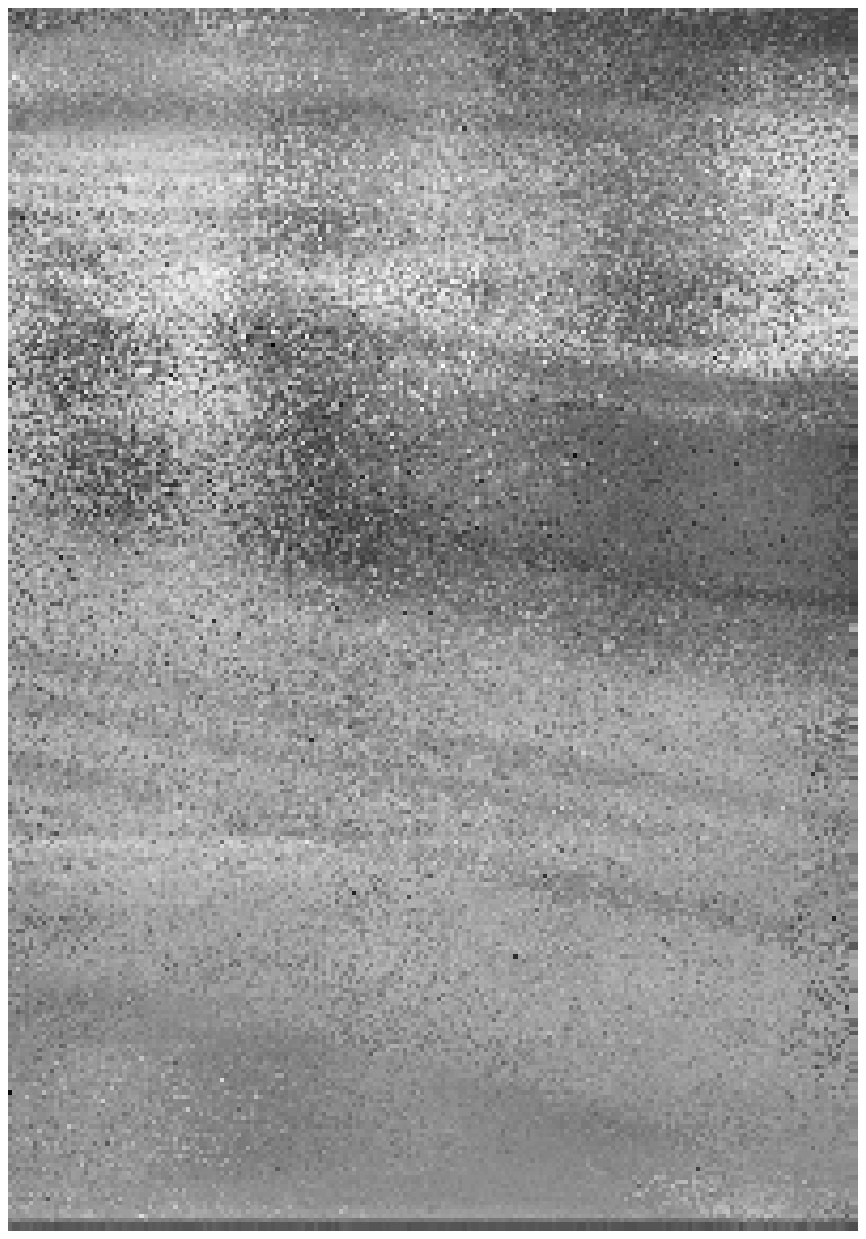}}
\vspace*{-2cm}
\centerline{\footnotesize (iv) PETRELS}\medskip
\end{minipage}
&%
\begin{minipage}[b]{.32\hsize}
\centerline{\includegraphics[width=1\hsize]{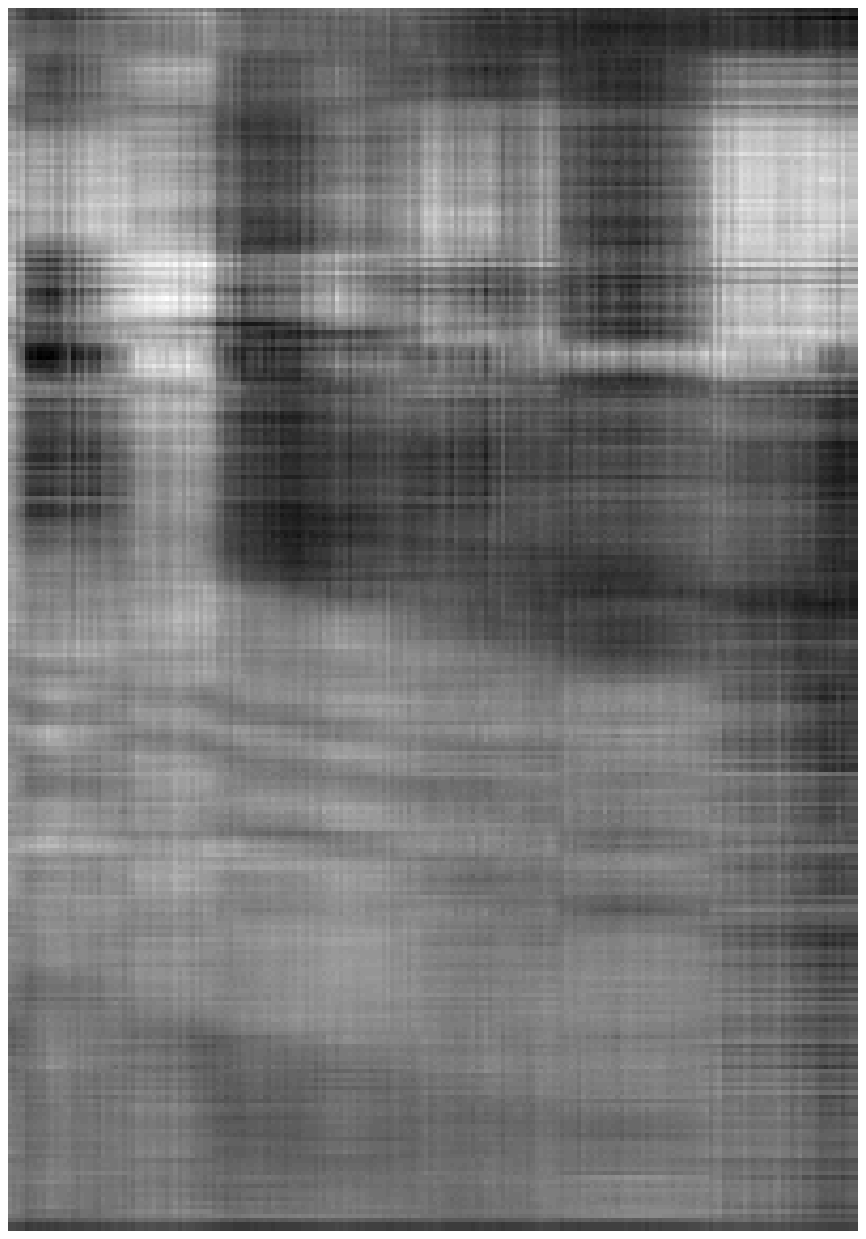}}
\vspace*{-2cm}
\centerline{\footnotesize (v) TeCPSGD}\medskip
\end{minipage}
&%
\begin{minipage}[b]{0.32\hsize}
\centerline{\includegraphics[width=1\hsize]{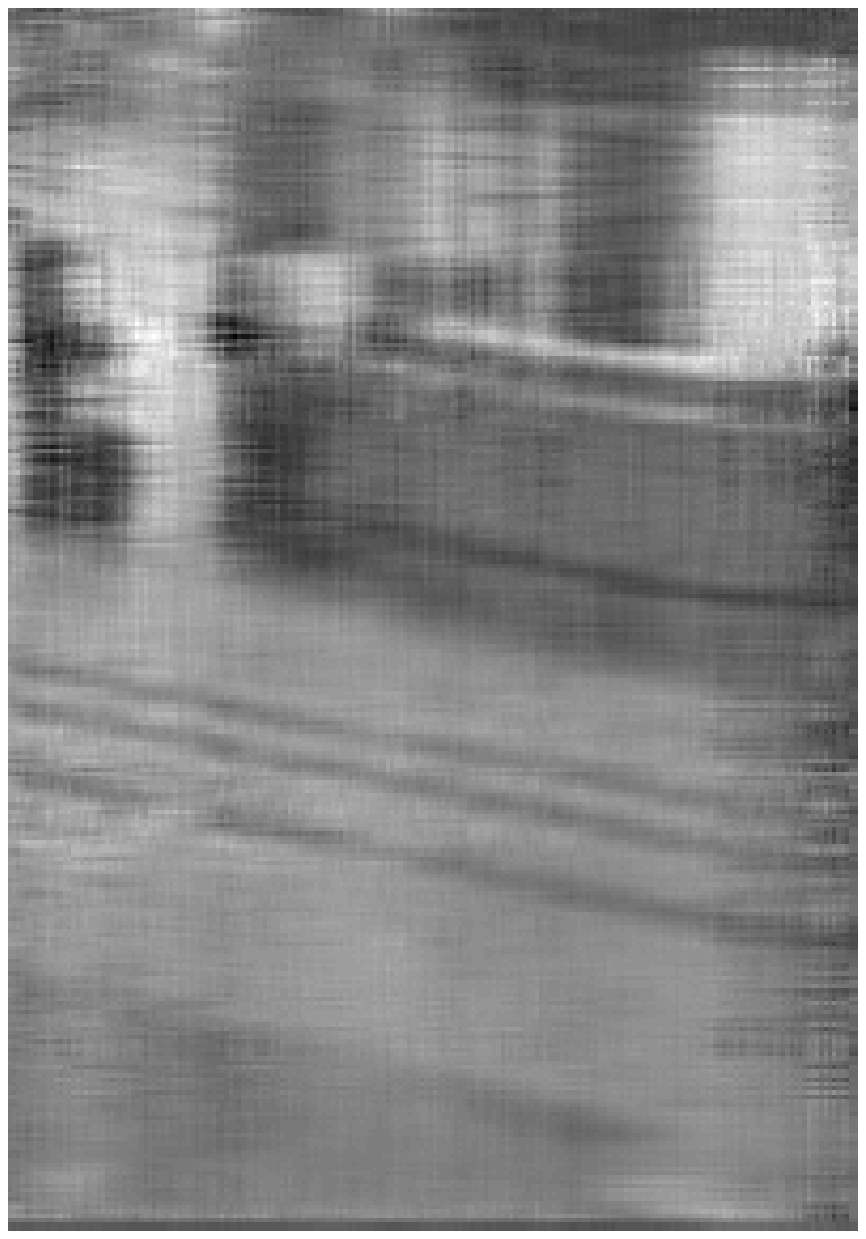}}
\vspace*{-2cm}
\centerline{\footnotesize (vi) OLSTEC}\medskip
\end{minipage}  
\end{tabular}
\centerline{(c) Reconstructed subspace images.}\medskip
\caption{The normalized estimation error in surveillance video Airport Hall.}
\label{fig:recimages}
\end{figure}

Finally, we evaluate the tracking performances using surveillance video. Although each video frame has no low-rank structure and a tensor-based approach basically presents shortcomings for the approximation of its underlying subspace, this section demonstrates the superior tacking performance of OLSTEC. We compare OLSTEC with TeCPSGD as well as the matrix-based algorithms including GROUSE \cite{Balzano_arXiv_2010_s}, GRASTA \cite{He_CVPR_2012}, and PETRELS \cite{Chi_IEEETransSP_2013}. ``Airport Hall" dataset of size $288\times 352$ with 500 frames is used. Moreover, for fair comparison between tensor and matrix-based algorithms, the rank is set to 20 and 5 for the former, i.e., OLSTEC and TeCPSGD, and for the latter, respectively. Still, the tensor-based algorithms have far fewer free parameters than those of the matrix-based algorithms. The observation ration $\rho$ is set to $0.1$. $\lambda$ and $\mu_r$ are $0.7$ and $0.1$, respectively, in the proposed algorithm. Actually, $\lambda$, $\mu$ and the stepsize are set, respectively, to $0.001$, $0.1$ and $10$ for TeCPSGD. 
This experiment particularly considers two scenarios.
The first separates moving objects in the foreground with static background. Figure \ref{fig:NRS_Video} (a) shows the superior performance of OLSTEC against other algorithms. 
The second scenario examines, by contrast, the performances against a dynamic moving background. The input video is created virtually by moving a cropped partial image from its original entire frame image of ``Airport Hall" video. The cropping window with $288 \times 200$ moves from the leftmost partial image to the rightmost, and returns to the leftmost image after stopping for a certain period of time. 
The generated video includes right-panning video from 38-th to 113-th frame and from 342-th to 417-th frame, and left-panning video from 190-th to 265-th frame.
Figure \ref{fig:NRS_Video}(b) shows how OLSTEC can quickly adapt to the changed background.
Figure \ref{fig:recimages} shows that the reconstructed image at 110-th frame  of OLSTEC gives better quality than that of others.  

\section{Conclusion and future work}

We have proposed a new online tensor subspace tracking algorithm, designated as OLSTEC, for a partially observed high-dimensional data stream corrupted by noise. Especially, we addressed a second-order stochastic gradient descent based on the recursive least squares to achieve faster convergence of subspace tracking. Numerical comparisons suggest that our proposed algorithm has superior performances for synthetic as well as real-world datasets.  As a future research direction, we will investigate the mechanisms of the Tucker decomposition.

\vfill\pagebreak

\bibliographystyle{unsrt}
\bibliography{/Users/kasai/Dropbox/DOC/Research/bibtex/matrix_tensor_completion,/Users/kasai/Dropbox/DOC/Research/bibtex/nips_matrix_tensor_completion,/Users/kasai/Dropbox/DOC/Research/bibtex/stochastic_online_learning}

\end{document}